\documentclass[11pt, reqno]{amsart}
\usepackage{mathrsfs}
\usepackage[active]{srcltx}
\usepackage{mathrsfs}
\usepackage{longtable}

\usepackage[T1]{fontenc}
\usepackage[latin9]{inputenc}
\usepackage[a4paper]{geometry}
\usepackage{amsthm}
\usepackage{amssymb}
\usepackage{mathtools}
\usepackage{mathrsfs}
\usepackage{xcolor}
\usepackage{amsmath,stmaryrd,graphicx}
\usepackage{hyperref}
\usepackage{enumerate}
\usepackage{titlesec}
\titlelabel{\thetitle.\quad}

\usepackage{graphicx}%
\usepackage{color}
\usepackage{cite}
\usepackage{hyperref}
\hypersetup{
	colorlinks = true,%
	citecolor = blue,
	filecolor=red,%
	linkcolor = [rgb]{0.65,0.0,0.0},%
	anchorcolor = red,
	pagecolor = red,
	urlcolor= [rgb]{0.65,0.0,0.0}
	linktocpage=true,
	pdfpagelabels=true,
	bookmarksnumbered=true,
}

\numberwithin{equation}{section}

\newtheorem{theorem}{Theorem}[section]

\theoremstyle{definition}

\newtheorem{proposition}{Proposition}[section]
\begin{document}

	\title[Differential inclusions on Riemannian manifolds]
	{Elliptic differential inclusions on non-compact Riemannian manifolds}
	
	\date{}

	\author{Alexandru Krist\'aly \and  Ildik\'o I. Mezei \and K\'aroly Szil\'ak}

	\thanks{A. Krist\'aly and I.I. Mezei are supported by the UEFISCDI/CNCS grant PN-III-P4-ID-PCE2020-1001.}

	\address{A. Krist\'aly: Department of Economics, Babe\c s-Bolyai University, Str. Teodor Mihali 58-60, 400591 Cluj-Napoca,
		Romania \& Institute of Applied Mathematics, \'Obuda
		University, B\'ecsi \'ut 96/B, 1034
		Budapest, Hungary}
	
	\email{alexandru.kristaly@ubbcluj.ro; kristaly.alexandru@nik.uni-obuda.hu}

	\address{I.I. Mezei: Department of Mathematics and Computer Sciences, Babe\c s-Bolyai University, Str. Mihail Kog\u alniceanu 1, 400084 Cluj-Napoca,
		Romania}
	\email{ildiko.mezei@ubbcluj.ro}
	
		\address{K. Szil\'ak: Institute of Applied Mathematics, \'Obuda
		University, B\'ecsi \'ut 96/B, 1034
		Budapest, Hungary}
	
	\email{karoly.szilak@gmail.com}

	\begin{abstract} 
		\noindent We investigate a large class of elliptic differential inclusions on non-compact complete Riemannian manifolds which involves the Laplace-Beltrami operator and a Hardy-type singular term.\ Depending on the behavior of the \textit{nonlinear term} and on the \textit{curvature} of the Riemannian manifold, we guarantee non-existence and existence/multi\-plicity of solutions for the studied differential inclusion. The proofs are based on nonsmooth variational analysis as well as  isometric actions and fine eigenvalue properties on Riemannian manifolds. The results are also new in the smooth setting. 
	\end{abstract}
	
	
	
	
	
\dedicatory{Dedicated to Professor Siegfried Carl on the occasion of his 70th birthday}
	
	\maketitle
	
	
	\vspace{-0.8cm}
	\section{Introduction} 
	
	Various geometric/physical phenomena can be reduced to finding solutions for the problem 
	$$\mathcal Lu(x)=\alpha(x)f(u(x)),\ x\in \Omega,\eqno{(P)}$$
	where $\Omega$ is an open domain in an ambient metric measure space,  $\mathcal L$ is an elliptic-type operator, $\alpha:\Omega\to \mathbb R$ is a measurable potential, and $f: \mathbb R\to \mathbb R$ is a nonlinear function having certain regularity and growth properties. Such problems arise from the Yamabe problem on compact/non-compact Riemannian manifolds, the standing Schr\"odinger  equation in $\mathbb R^n$ $(n\geq 2)$, Dirichlet and Neumann problems on bounded/unbounded domains, etc. Wide range of strategies and theories have been applied in the last century in order to investigate problem $(P)$, as variational methods, fixed point arguments, sub- and super-solution techniques, etc.   
	
	An important class of problems within $(P)$ appears when the nonlinear term $f(x,\cdot)$ is not necessarily continuous; such a relevant example appears in the description of the von K\'arm\'an adhesive plates, see Panagiotopoulos \cite{Panagiotopoulos}. Due to the jumping effect of $f(x,\cdot)$,  as a first approach, problem $(P)$ need not has any solution. However, from physical  reasons, we expect to obtain certain equilibrium states of the phenomena described by means of problem $(P)$. Accordingly, a natural way to handle the aforementioned discontinuity situation  is to 'fill the gaps', defining a differential inclusion associated with problem $(P)$. More precisely, if $f$ is locally essentially bounded on $\mathbb R$, we consider instead of the value $f(t)$ the interval $[\underline f(t),\overline f(t)]$, where 
	$$\underline f(t)=\lim_{\delta\to 0^+}{\rm essinf}_{|s-t|<\delta}f(s),\ \ \overline f(t)=\lim_{\delta\to 0^+}{\rm esssup}_{|s-t|<\delta}f(s);$$
	here, ${\rm essinf}_A f=\sup\{a\in \mathbb R:f(x)\geq a\ {\rm for\ a.e.}\ x\in A\}$ and ${\rm esssup}_A f=-{\rm essinf}_A (-f)$ whenever  $A\neq \emptyset$. In this way, we replace $(P)$ by the differential inclusion problem 
	$$\mathcal Lu(x)\in \alpha(x)\partial F(u(x)),\ x\in \Omega,\eqno{(DI)}$$
	where  $F(t)=\displaystyle\int_0^t f(s)ds$ is a locally Lipschitz function, and $\partial F(t)=[\underline f(t),\overline f(t)]$, $t\in \mathbb R$. Hereafter, $\partial F$ stands for the  subdifferential of $F$ at $t\in \mathbb R$ in the sense of Clarke \cite{Clarke}.
	
	Differential inclusion problems, similar to $(DI),$ may appear on not necessarily Euclidean structures; indeed, in certain circumstances the domain $\Omega$ can be a subset of a \textit{curved space} (Riemannian or Finsler manifolds,  sub-Riemannian structures as Heisenberg or Carnot groups, etc.), while the operator $\mathcal L$ may reflect the geometric feature of the ambient space.  
	
	In the present paper we consider the variational differential inclusion  
	\begin{equation}\label{problem-1}
	\mathcal L_{x_0}u(x)\equiv	-\Delta_gu(x) -\mu\frac{u(x)}{d^2_g(x_0,x)}+u(x)\in \lambda\alpha(x)\partial F(u(x)), \ \ x\in M,
	\end{equation}
	where $(M,g)$ is an $n$-dimensional complete Riemannian manifold, $n\geq 3$ (endowed with its canonical measure ${\rm d}v_g$), $\Delta_g$ is the Laplace-Beltrami operator on $(M,g)$, $d_g:M\times M\to \mathbb R$ is the distance function associated with the Riemannian metric $g$,  $x_0\in M$ is a fixed point, $\alpha:M\to \mathbb R$ is a measurable potential, $\mu,\lambda\in \mathbb R$ are some parameters,  $F:\mathbb R\to \mathbb R$ is a locally Lipschitz function and $\partial F$ stands for the Clarke subdifferential of $F$. An element $u\in H^1(M)$  is a solution of \eqref{problem-1} if there exists a measurable selection $x\mapsto \xi_x\in \partial F(u(x))$ such that  the map $x\mapsto \alpha(x)\xi_x w(x)$ belongs to $L^1(M)$ for every test-function $w\in H^1(M)$ and one has
	\begin{equation}\label{solution-DI}
	\int_M \nabla_gu(x)\nabla_gw(x){\rm d}v_g -\mu\int_M\frac{u(x)w(x)}{d^2_g(x_0,x)}{\rm d}v_g+\int_Mu(x)w(x){\rm d}v_g=\lambda\int_M \alpha(x)\xi_x w(x){\rm d}v_g.
	\end{equation}
One can readily observe that \eqref{solution-DI} reduces to  the fact that $u$ is a weak solution of 
	\begin{equation}\label{problem-2}
	-\Delta_gu(x) -\mu\frac{u(x)}{d^2_g(x_0,x)}+u(x)= \lambda\alpha(x)f(u(x)), \ \ x\in M,
\end{equation}
 whenever $f$ is continuous (and consequently, $F$ is of class $C^1$ and  $\partial F(t)=F'(t)=f(t)$). 
 
On one hand, variational elliptic differential inclusions  as \eqref{problem-1} -- or  slightly different versions of them formulated in terms of variational-hemivariational inequalities  --  have been deeply studied in the last three decades, mostly in Euclidean spaces (both for bounded and unbounded domains), see e.g.\  Bonanno, D'Agu\`i and Winkert \cite{Bonanno-etal},  Candito and  Livrea \cite{Candito-Livrea}, Carl and Le \cite{Carl-Le, Carl-Le-2}, Carl,  Le and Motreanu \cite{CLM},  Costea, Krist\'aly and Varga \cite{CKV}, Gasi\'nski and  Papageorgiou \cite{Gasinski-Papa}, Krist\'aly and Varga \cite{KV}, Liu,  Liu and Motreanu \cite{Liu-Liu-Motreanu}, 
Liu,  Livrea,  Motreanu and Zeng \cite{olaszok-1}, Mig\'orski, Ochal and Sofonea \cite{MOS}, Motreanu and Panagiotopoulos \cite{Motreanu-Pan}, Panagiotopoulos \cite{Panagiotopoulos}, Varga \cite{Varga}, etc.\ On the other hand, various forms of \eqref{problem-2} have been investigated  both on compact and non-compact Riemannian manifolds (mostly without the singular term), see e.g.\ Berchio,  Ferrero and Grillo \cite{Berchio}, Bonanno, Molica Bisci and R\u adulescu \cite{Bonanno-Bisci}, Jaber \cite{Jaber}, Lisei and  Varga \cite{Hanne-Varga}, Liu and  Liu \cite{Liu-Liu}, 	   Molica Bisci and Pucci \cite{BPucci}, Molica Bisci and Secchi \cite{Molica-0}, Molica Bisci and Vilasi \cite{Molica}, etc.\ As expected, on non-compact manifolds additional restrictions and approaches are needed  to compensate the lack of compactness.  

We shall focus to a broad class of \textit{non-compact Riemannian manifolds} and  prove various non-existence, existence and multiplicity results for the differential inclusion problem \eqref{problem-1}, by assuming certain  \textit{curvature} hypotheses and \textit{growths} for the function  $F$  (at the origin and at infinity).  In fact, we  consider two  classes of Riemannian manifolds having different curvature restrictions; namely, we assume that a complete, non-compact Riemannian manifold $(M,g)$ satisfies one of the conditions:   
	\begin{itemize}
		\item[(i)]  ${\bf K}\leq -\kappa$ for some $\kappa\geq 0$, where ${\bf K}$ is the sectional curvature of the Cartan-Hadamard manifold\footnote{Complete, simply connected Riemannian manifold with non-positive sectional curvature.} $(M,g)$;
		\item[(ii)]  {\sf Ric}$_{(M,g)}\geq 0$, where  {\sf Ric}$_{(M,g)}$ is the Ricci curvature on  $(M,g)$.
	\end{itemize}
The 'clash' of (i) and (ii) is precisely the Euclidean space $\mathbb R^n$ endowed with the usual metric. In the case (ii), i.e., when {\sf Ric}$_{(M,g)}\geq 0$, a crucial role is played by  the  \textit{asymptotic volume ratio} 
$${\sf AVR}_{(M,g)}=\lim_{r\to \infty}\frac{{V}_g(B_x(r))}{\omega_nr^n},$$
where $V_g$ stands for the volume in $(M,g)$, $B_x(r)=\{y\in M:d_g(x,y)<r\}$ is the ball of radius $r>0$ and center $x\in M$, while $\omega_n=\pi^{n/2}/\Gamma(1+n/2)$ is the volume of the Euclidean unit ball in $\mathbb R^n$. By Bishop-Gromov comparison principle it turns out that the asymptotic volume ratio is well-defined (i.e.,  independent of the choice of $x\in M$) and ${\sf AVR}_{(M,g)}\in [0,1]$.

	We assume on the potential $\alpha:M\to \mathbb R$ that 
	
	$(\textbf{H})_\alpha:$ $\alpha\geq 0$ and $\alpha\in L^1(M)\cap L^\infty(M)\setminus \{0\}$.  
	
\noindent For the locally Lipschitz function $F:\mathbb R\to \mathbb R$ we require 

	$(\textbf{H})_0:$ there exists $C_0>0$ such that $$|\xi|\leq C_0 |t|,\ \ \forall \xi\in \partial F(t),\ t\in \mathbb R.$$
	
The first result of the paper reads as follows.

\begin{theorem}\label{theorem-non-existence}{\rm (Non-existence)} Let $(M,g)$ be an $n$-dimensional complete non-compact Riemannian manifold, $n\geq 3$, and assume that the potential $\alpha:M\to \mathbb R$ and the locally Lipschitz function $F:\mathbb R\to \mathbb R$ satisfy assumptions	$( \bf{H})_\alpha$ and    $({\bf H})_0$, respectively. Assume in addition that one of the following curvature conditions holds$:$
	\begin{itemize}
		\item[(i)] ${\bf K}\leq -\kappa$ for some $\kappa\geq 0,$ $(M,g)$ is simply connected and 
		\begin{itemize}
			\item[(i1)] either $\kappa=0$,  $\mu\leq \frac{(n-2)^2}{4}$ and $|\lambda|C_0\|\alpha\|_{L^\infty}\leq 1,$ 
			\item[(i2)] or $\kappa>0$, $\mu\leq \frac{(n-2)^2}{4}$ and $(n-2)^2(|\lambda|C_0\|\alpha\|_{L^\infty}-1)\leq (n-1)^2\left(\frac{(n-2)^2}{4}-\mu_+\right)\kappa,$ where $\mu_+=\max(\mu,0);$
		\end{itemize}
	\item[(ii)] {\rm {\sf Ric}}$_{(M,g)}\geq 0$, $\mu\leq {\sf AVR}_{(M,g)}^\frac{2}{n}\frac{(n-2)^2}{4}$ and $|\lambda|C_0\|\alpha\|_{L^\infty}\leq 1.$ 
	\end{itemize}
Then the differential inclusion \eqref{problem-1} has only the zero solution.
\end{theorem}
	
	The assertions in (i) show that there is a balance in the sense that when a stronger curvature restriction occurs, the analytic assumption can be relaxed. The proof of Theorem \ref{theorem-non-existence} is based on a direct computation combined with Hardy-type inequalities and sharp spectral gap  estimates on Riemannian manifolds; the ingredients to the proof are recalled in \S \ref{section-2-1}.  
	
In order to produce existence or even multiplicity of non-zero solutions to \eqref{problem-1}, we require on the locally Lipschitz function $F:\mathbb R\to \mathbb R$ the following assumptions: 

		$(\textbf{H})_1:$ 
		 $\displaystyle\lim_{t\to 0}\frac{\max\{|\xi|:\xi\in \partial F(t)\}}{t}=0;$
		 
		 	$(\textbf{H})_2:$ 
		 $\displaystyle\lim_{|t|\to \infty}\frac{\max\{|\xi|:\xi\in \partial F(t)\}}{t}=0;$
		 
	
		$(\textbf{H})_3:$ $F(0)=0$ and there exist $t_0^-<0<t_0^+$ such that 
	$F(t_0^\pm)>0$.
	
\noindent 	Note that $(\textbf{H})_1$ and $(\textbf{H})_2$ mean that the function $t\mapsto \max\{|\xi|:\xi\in \partial F(t)\}$ is \textit{superlinear at the origin} and \textit{sublinear at infinity}, respectively; in particular, by using Lebourg's mean value theorem, we observe that  $F$ is \textit{sub-quadratic at infinity}.  
In addition, by the upper semicontinuity of the set-valued function $t\mapsto \partial F(t)$ and conditions $(\textbf{H})_1$ and $(\textbf{H})_2$, it turns out that the hypothesis $(\textbf{H})_0$ is also valid for a suitably large value of $C_0>0$; in particular, Theorem \ref{theorem-non-existence} can be applied (under the assumptions $(\textbf{H})_1$ and $(\textbf{H})_2$), and  for sufficiently 'small' values of $|\lambda|$ only the zero solution exists for the differential inclusion \eqref{problem-1}. 
However, for 'large' values of $\lambda>0$, we can guarantee the existence of multiple non-zero solutions for \eqref{problem-1} by requiring further assumptions on the behavior of the \textit{isometric group} of the Riemannian manifold $(M,g)$. In fact, the latter assumptions are destined to balance the lack of compactness of the Riemannian manifolds we are dealing with. 

To state the second result of the paper, we denote by ${\sf Isom}_g(M )$ the\textit{ group of isometries} of the complete Riemannian manifold $(M,g)$. Let   $G$ be a subgroup of ${\sf Isom}_g(M )$ and
\begin{equation}\label{fixed-points}
	{\sf Fix}_M(G)=\{x\in M:\sigma (x)=x,\forall \sigma\in G\}
\end{equation}
be the set of \textit{fixed points} of the isometry group  $G$ in $M$.  The $G$-\textit{orbit} of a point $x\in M$  is $\mathcal O_G^x=\{\sigma(x):\sigma\in G\}$. The continuous action of the group $G$ on $M$ is \textit{coercive} if for every $t > 0$ the set $\mathcal O_t := \{x \in M:{\rm diam} (\mathcal O^x_G) \leq t\}$ is bounded, see Skrzypczak and Tintarev \cite{ST-1,ST-2}; here diam$(S)$ denotes the diameter of $S\subset M$. A function $u:M\to \mathbb R$ is $G$-\textit{invariant} if $u(x)=u(\sigma(x))$ for every $x\in M$ and $\sigma\in G.$ 

\begin{theorem}\label{theorem-multiplicity}{\rm (Multiplicity:\ sub-quadratic nonlinearity at infinity)} Let $(M,g)$ be an $n$-dimen\-sional complete non-compact Riemannian manifold, $n\geq 3$, and $G$ be a compact connected subgroup of ${\sf Isom}_g(M )$ such that ${\sf Fix}_M(G) =\{x_0\}$ for the same $x_0\in M$ as in problem \eqref{problem-1}. Let $\alpha:M\to \mathbb R$ be a  potential satisfying $(\bf{H})_\alpha$ which depends only on $d_g(x_0,\cdot)$ and the locally Lipschitz function $F:\mathbb R\to \mathbb R$ satisfying assumptions	   $({\bf H})_i$, $i\in\{1,2,3\},$ respectively. In addition, we assume that one of the following curvature assumptions holds$:$
	\begin{itemize}
		\item[(i)] $(M,g)$ is of Cartan-Hadamard-type and 
		$0\leq \mu< \frac{(n-2)^2}{4};$ 
		\item[(ii)] {\rm {\sf Ric}}$_{(M,g)}\geq 0,$ ${\sf AVR}_{(M,g)}>0,$ $0\leq \mu< {\sf AVR}_{(M,g)}^\frac{2}{n}\frac{(n-2)^2}{4}$  and $G$ is coercive. 
	\end{itemize}
	Then there exists $\lambda_0>0$ such that for every $\lambda>\lambda_0$ the differential inclusion \eqref{problem-1} has at least four non-zero $G$-invariant solutions in $H^1(M)$.
\end{theorem}

The proof of Theorem \ref{theorem-multiplicity} is based on truncation and variational arguments, combined with careful isometry actions on $H^1(M)$. The key ingredients are the nonsmooth principle of symmetric criticality and mountain pass theorem (with the Palais-Smale  condition) and the compact embedding of $G$-invariant functions of $H^1(M)$ into appropriate Lebesgue spaces over $M$, which are valid in both geometric contexts (i) and (ii).  Examples of Riemannian manifolds  with the above curvature restrictions  and isometric actions are presented  in Krist\'aly \cite{Kristaly-JFA} and  Farkas, Krist\'aly and Mester \cite{FKM} in the setting (i), and Balogh and Krist\'aly \cite{Balogh-Kristaly-Annalen} in the  framework (ii).

As we already noticed, assumptions $(\textbf{H})_1$ and $(\textbf{H})_2$ imply that  $F$ is sub-quadratic at infinity. In the sequel, we establish a counterpart of Theorem \ref{theorem-multiplicity} whenever $F$ is \textit{super-quadratic at infinity}. More precisely, we assume that the locally Lipschitz function $F:\mathbb R\to \mathbb R$ satisfies the following assumptions: 

	$(\textbf{H})_4:$ $F(0)=0$ and there exist  $\nu>2$ and $C>0$  such that  
	\begin{equation}\label{cerami}
			2F(t)+F^0(t;-t)\leq -C|t|^{\nu},\ \ \forall t\in \mathbb R;
	\end{equation}

$(\textbf{H})_5:$ there is $q\in \left(2,2+\frac{4}{n}\right)$ such that $\max\{|\xi|:\xi\in \partial F(t)\}=O(|t|^{q-1})$ as $|t|\to \infty$.

\noindent Here, $F^0(t;s)$ is the generalized directional derivative of $F$ at the point $t\in \mathbb R$ and direction $s\in \mathbb R$, see \S \ref{section-2}. Note that by $(\textbf{H})_1$ and $(\textbf{H})_4$, $F$
 is super-quadratic at infinity, see \S \ref{section-4}.   
  
 \begin{theorem}\label{theorem-multiplicity-2}{\rm (Existence/Multiplicity:\ super-quadratic nonlinearity at infinity)} Let $(M,g)$ be an $n$-dimensional complete non-compact Riemannian manifold, $n\geq 3$, and $G$ be a compact connected subgroup of ${\sf Isom}_g(M )$ such that ${\sf Fix}_M(G) =\{x_0\}$ for the same $x_0\in M$ as in problem \eqref{problem-1}. Let $\alpha\in L^\infty(M)$ be a  potential which depends only on $d_g(x_0,\cdot)$ and ${\rm essinf}_{x\in M}\alpha(x)=\alpha_0>0$, while the locally Lipschitz function $F:\mathbb R\to \mathbb R$ satisfies the  assumptions	   $({\bf H})_i$, $i\in\{1,4,5\},$ respectively. If one of the curvature assumptions {\rm (i)} or {\rm (ii)} holds from Theorem \ref{theorem-multiplicity}, 
%
 	then for every $\lambda>0$ the differential inclusion \eqref{problem-1} has at least a non-zero $G$-invariant solution in $H^1(M)$. In addition, if $F$ is an even function, \eqref{problem-1} has infinitely many distinct $G$-invariant solutions in $H^1(M)$. 
 \end{theorem}
 
The proof is based on the same geometric arguments as in Theorem \ref{theorem-multiplicity} (curvature constraints, isometric actions), combined with the nonsmooth mountain pass  or fountain theorem involving the Cerami compactness condition. 

The paper is organized as follows. In \S \ref{section-prelim} we collect  those results that are indispensable in our proofs. Namely, we first recall certain functional inequalities and spectral estimates on Riemannian manifolds; then we recall some elements from the nonsmooth calculus of Clarke \cite{Clarke} including also the nonsmooth principle of symmetric criticality. In \S \ref{section-3} we prove the non-existence results, established within Theorem \ref{theorem-non-existence}. In \S \ref{section-4} we discuss our first existence/multiplicity results in the  sub-quadratic case, by proving  Theorem \ref{theorem-multiplicity}. Finally, Section \S \ref{section-5} is devoted to handle the  super-quadratic case, i.e.,  Theorem \ref{theorem-multiplicity-2}. 

\section{Preliminaries}\label{section-prelim} 
In this section we recall those notions and results that are crucial to carry out our proofs. Before to do this, we fix some notations. If $(M,g)$ is a complete Riemannian manifold,  the Sobolev space
 $H^1(M)$ over $M$ is the completion of $C_0^\infty(M)$ with respect to the norm 
 $$\|u\|_{H_1}=\left(\int_M |\nabla_g u|^2{\rm d}v_g+\int_Mu^2{\rm d}v_g\right)^{1/2},$$
 while the $L^q$-Lebesgue norm $(q\geq 1)$ is 
 $$\|u\|_{L^q}=\left(\int_M|u|^q{\rm d}v_g\right)^{1/q},$$
 with the  supremum-norm for $q=+\infty.$
 
\subsection{Functional inequalities and spectral estimates on  Riemannian manifolds.}\label{section-2-1} 
\subsubsection{Cartan-Hadamard manifolds.}\label{subsub-2-1-1} Throughout this subsection, let $(M,g)$ be an $n$-dimensional Cartan-Hadamard manifold, $n\geq 3$.  We notice that in this geometric context, there exists $C_n>0$  such that
$$\|u\|_{L^{2^*}}\leq C_{n}\left(\int_M |\nabla_g u|^2{\rm d}v_g\right)^{1/2},\ \ \forall u\in C_0^\infty(M),$$
see e.g. Hebey \cite[Chapter 8]{Hebey}, where $2^*=2n/(n-2)$ is the  critical Sobolev exponent. Moreover, the best Sobolev embedding constant $C_n$ is precisely its Euclidean counterpart ${\sf AT}_n$, provided by Aubin \cite{Aubin} and Talenti \cite{Talenti},  whenever the Cartan-Hadamard conjecture holds on $(M,g)$ (e.g. in dimensions 3 and 4). In high-dimensions, the sharp constant $C_n>0$ is not known; however, a non-optimal form can be given by means of the Croke-constant as in Hebey \cite[p.\ 239]{Hebey}. 

A density argument combined with a simple interpolation shows that the Sobolev space $H^1(M)$  is continuously embedded into $L^q(M)$ for every $q\in [2,2^*]$; more precisely, there exists $K_q^->0$ such that
\begin{equation}\label{Sobolev-CH}
	\|u\|_{L^q}\leq K_q^-\|u\|_{H^1},\ \ \forall u\in H^1(M).
\end{equation} 

Let $x_0\in M$ be fixed. Then the \textit{Hardy inequality} holds on $(M,g)$, which reads as
\begin{equation}\label{Hardy-CH}
\frac{(n-2)^2}{4}	\int_M \frac{u^2(x)}{d_g^2(x_0,x)}{\rm d}v_g\leq \int_M |\nabla_g u|^2{\rm d}v_g,\ \ \forall u\in H^1(M),
\end{equation}
where $\frac{(n-2)^2}{4}$ is sharp and never achieved, see e.g. D'Ambrosio and Dipierro \cite{Dipierro}, and Krist\'aly \cite{Kristaly-JMPA}. 

In addition, if  the sectional curvature has the property ${\bf K}\leq -\kappa$ for some $\kappa> 0,$ then  \textit{McKean's  spectral gap theorem} asserts that
\begin{equation}\label{McKean-spectral-gap}
	\gamma_{(M,g)}:=\inf_{u\in H^1(M)\setminus \{0\}}\frac{\displaystyle\int_M |\nabla_g u|^2{\rm d}v_g}{\displaystyle\int_M u^2{\rm d}v_g}\geq \frac{(n-1)^2}{4}\kappa.
\end{equation}
The inequality  \eqref{McKean-spectral-gap} is sharp, see e.g.\ on the $n$-dimensional hyperbolic space $\mathbb H^n_\kappa$ with constant sectional curvature ${\bf K}= -\kappa$; we also notice that the infimum in \eqref{McKean-spectral-gap}
is not achieved by any function $u\in H^1(M)$; 

\subsubsection{Riemannian manifolds with non-negative Ricci curvature.} In this subsection we consider an $n$-dimensional $(n\geq 3)$ complete non-compact Riemannian manifold $(M,g)$ with ${\rm {\sf Ric}}$$_{(M,g)}\geq 0$. As we already noticed in the Introduction, the asymptotic volume ratio ${\sf AVR}_{(M,g)}\in [0,1]$ provides deep geometric information about the manifold; for instance,   ${\sf AVR}_{(M,g)}=1$ if and only if $(M,g)$ is isometric to the Euclidean space $\mathbb R^n$. Quantitatively speaking, closer value of ${\sf AVR}_{(M,g)}$ to 1 implies  topologically closer manifold $(M,g)$ to the Euclidean space $\mathbb R^n$, expressed in terms of the trivialization of higher homotopy groups of $M$, see Munn \cite{Munn}.  

In the geometric context when $(M,g)$ is a complete non-compact Riemannian manifold with {\rm {\sf Ric}}$_{(M,g)}\geq 0$, a necessarily and sufficient condition to have the Sobolev embedding is the fact that  ${\sf AVR}_{(M,g)}>0$, see 
Coulhon and Saloff-Coste \cite{Coulhon-SC}
 and Hebey \cite{Hebey}. Moreover, a recent result of Balogh and Krist\'aly \cite{Balogh-Kristaly-Annalen} asserts that if ${\sf AVR}_{(M,g)}>0$ then 
$$\|u\|_{L^{2^*}}\leq {\sf AVR}_{(M,g)}^{-\frac{1}{n}}{\sf AT}_n\left(\int_M |\nabla_g u|^2{\rm d}v_g\right)^{1/2},\ \ \forall u\in H^1(M),$$
where the constant ${\sf AVR}_{(M,g)}^{-\frac{1}{n}}{\sf AT}_n$ is sharp; here, as in \S \ref{subsub-2-1-1},   ${\sf AT}_n$ stands for the best Sobolev embedding constant in the Euclidean Sobolev inequality on $\mathbb R^n$, see Aubin \cite{Aubin} and Talenti \cite{Talenti}. In particular, $H^1(M)$  is continuously embedded into $L^q(M)$ for every $q\in [2,2^*]$; more precisely, there exists $K_q^+>0$ such that
\begin{equation}\label{Sobolev-Ric}
	\|u\|_{L^q}\leq K^+_q\|u\|_{H^1},\ \ \forall u\in H^1(M).
\end{equation}  

Given $x_0\in M$ fixed, the \textit{Hardy inequality}  on $(M,g)$ is verified as 
\begin{equation}\label{Hardy-Ric}
	{\sf AVR}_{(M,g)}^\frac{2}{n}\frac{(n-2)^2}{4}	\int_M \frac{u^2(x)}{d_g^2(x_0,x)}{\rm d}v_g\leq \int_M |\nabla_g u|^2{\rm d}v_g,\ \ \forall u\in H^1(M), 
\end{equation}
see Krist\'aly, Mester and Mezei  \cite{KMM}. The sharpness of the constant in \eqref{Hardy-Ric} is not known unless we are in the classical Euclidean setting. However, if we  assume that there exists a non-zero function realizing the equality in \eqref{Hardy-Ric}, from the proof in \cite{KMM}  (based on a P\'olya-Szeg\H o inequality involving the asymptotic volume ratio) we would obtain that ${\sf AVR}_{(M,g)}=1$, i.e., $(M,g)$ is isometric to the Euclidean space $\mathbb R^n$, which is a contradiction (since no non-zero extremal function exists in the Hardy inequality in $\mathbb R^n$).

\subsubsection{Compact embeddings via isometric actions.}\label{subsection-compact} According to the previous subsections, the Sobolev space $H^1(M)$ is continuously embedded into $L^q(M)$, $q\in [2,2^*]$, whenever $(M,g)$ is either a Cartan-Hadamard manifold or a complete non-compact Riemannian manifold with ${\rm {\sf Ric}}$$_{(M,g)}\geq 0$; however, none of them   is compact, which represents an impediment to apply variational arguments on $H^1(M)$. To handle the lack of compactness, we use certain symmetrization \`a la Lions \cite{Lions} by means of \textit{isometries} of $M$. 

Let $(M,g)$ be an $n$-dimensional complete non-compact Riemanian manifold, $n\geq 3$, and -- as in the Introduction -- ${\sf Isom}_g(M )$ be the group of isometries of $(M,g)$. Let  $G$ be a subgroup of ${\sf Isom}_g(M )$ and 	${\sf Fix}_M(G)$ be the set of fixed points of the isometry group  $G$ in $M$, see \eqref{fixed-points}. 
Let $$H_G^1(M)=\{u\in H^1(M): u\circ \sigma =u,\ \forall \sigma\in G\}$$
be the closed $G$-invariant subspace of $H^1(M).$ The consequences of the main results in the paper by Farkas, Krist\'aly and Mester \cite{FKM} state that if one of the following assumptions hold, i.e.,
\begin{itemize}
	\item $(M,g)$ is a Cartan-Hadamard manifold and ${\sf Fix}_M(G)$ is a singleton, or
	\item  {\rm {\sf Ric}}$_{(M,g)}\geq 0,$ ${\sf AVR}_{(M,g)}>0$ and $G$ is coercive,
\end{itemize} 
then the space $H_G^1(M)$ can be compactly embedded into $L^q(M)$ for every $q\in (2,2^*)$. 

\subsection{Non-smooth analysis.}\label{section-2}
\subsubsection{Locally Lipschitz functions.} In this subsection we recall those basic properties of locally Lipschitz functions which are used in our proofs; for details, see Clarke \cite{Clarke}. 

Let $X$ be a real Banach space with the norm $\|\cdot\|$.  A function $h:X \rightarrow\mathbb{R}$ is  {\it locally Lipschitz} if every point  $u\in X$ possesses a neighborhood  $U_{u} \subset X$ such that
\begin{equation}\label{loc-lip}
	\vert
	h(u_{1})-h(u_{2})\vert\leq K\| u_{1}-u_{2}\|, \quad \forall u_{1}, u_{2} \in U_{u},
\end{equation}
for a constant $K > 0$ depending on $U_{u}$. The {\it generalized directional derivative} of the locally
Lipschitz function $h:X\rightarrow\mathbb{R}$ at  $u\in
X$ in the direction $v\in X$ is given by
$$
h^{0}(u;v):=\limsup_{\scriptstyle {\it w\rightarrow u}\atop
	\scriptstyle \it t\searrow 0}\frac{h(w+tv)-h(w)}{t}.
$$ If  $h:X\rightarrow\mathbb{R}$ is a function of class $C^1$ on $X$, then
$	h^{0}(u;v)=\langle h'(u), v\rangle$ for all $u, v\in X.$ Hereafter, $\langle\cdot,\cdot\rangle$ and $\|\cdot\|_*$ stand for the duality mapping on $(X^*,X)$ and the norm on $X^*$, respectively. 
The {\it Clarke subdifferential\index{Clarke!subdifferential\}} $\partial h(u)$ of $h$ at a
	point $u\in X$  is the subset of the dual space $X^{\ast}$ given
	by 
	$$
	\partial h(u):=\left\{\zeta\in X^{\ast}: \langle \zeta, v\rangle
	\leq h^{0}(u; v),\ \forall v\in X \right\}.
	$$
	An element	$u\in X$ is a \textit{critical point} of $h$ if $0\in \partial h(u)$, see Chang \cite[Definition 2.1]{Chang}.  
	
	\begin{proposition}\label{prop-lok-lip-0} (Clarke \cite{Clarke}) \it 
		Let $h:X\rightarrow\mathbb{R}$ be a locally Lipschitz function. The following assertions hold$:$
		\begin{enumerate} [$(i)$]
			\item For every $u\in X$, $\partial h(u)$ is a nonempty, convex
			and weak$^{\ast}$-compact subset of $X^{\ast}$. Moreover,
			$
			\|\zeta\|_\ast\leq K$ for all $ \zeta\in \partial h(u),
			$
			with $K>0$ from \eqref{loc-lip}.
			
			\item For every $u\in X,$ $h^{0}(u; \cdot)$ is the support
			function of $\partial h(u)$, i.e.,
			$$
			h^{0}(u; v)=\max\left\{\langle \zeta, v\rangle:  \zeta\in \partial h(u)\right\},\ \ 
			\forall v\in X.
			$$
			
			\item The set-valued map
			$\partial h:X\rightsquigarrow {X^{\ast}}$ is closed from $s-X$ into $w^\ast-X^\ast$.    
			 In particular, if $X$ is finite dimensional, then $\partial h$ is an upper semicontinuous set-valued map.
			 
			 \item $($Lebourg's mean value theorem$)$ Let  $U$  be an open subset of a Banach space $X$ and $u, v$ be
			 two points of $U$ such that the line segment  $ [u, v]=\{
			 (1-t)u+tv : 0\leq t\leq 1 \}\subset U$. If  $h:U\rightarrow\mathbb{R}$ is a Lipschitz function, then
			 there exist $w\in (u,v)$ and $\zeta\in \partial h(w)$ such that 
			 $
			 h(v)-h(u)=\langle \zeta,v-u\rangle.
			 $ 
			 
			 \item If $j:X\to \mathbb R$ is of class $C^1$ on $X$, then $\partial (j+h)(u)=j'(u)+\partial h(u)$ and 
			 $(j+h)^0(u;v)=\langle j'(u),v\rangle +h^0(u;v)$ for every $u,v\in X.$ 
			 \item $(-h)^0(u;v)=h^0(u;-v)$ for every $u,v\in X;$ 
			 \item $\partial (sh)(u)=s\partial h(u)$ for every $s\in \mathbb R$ and $u\in X.$
			
		\end{enumerate}
	\end{proposition}


\subsubsection{Principle of symmetric criticality for locally Lipschitz functionals.}

Let $G$ be a compact Lie group acting \textit{linear isometrically} on the real Banach
space $(X, \|\cdot\|)$, i.e., the action $G\times X\to X$, $(\sigma ,u)\mapsto \sigma u$ is continuous and
for every $\sigma \in G$ the map $u\mapsto \sigma u$ is linear such that $\|\sigma u\|=\|u\|$ for every $u\in X$.
A function $h:X\to \mathbb R$ is $G$-\textit{invariant} if $h(\sigma u) = h(u)$ for all $\sigma \in G$ and $u\in X$.
Let 
\begin{equation}\label{X-fixed-point}
	{\sf Fix}_X(G) = \{u \in  X : \sigma u = u, \ \forall \sigma  \in  G\}
\end{equation}
be the \textit{set of fixed points} of $G$ over $X$. According to Krawcewicz and
Marzantowicz \cite{KM} (see also Costea, Krist\'aly and Varga \cite[Section 3.4]{CKV}),   the \textit{principle of
symmetric criticality} for locally Lipschitz functions can be stated as follows. 

\begin{proposition}\label{psc} (Krawcewicz and
	Marzantowicz \cite{KM})  \it Let $G$ be a compact Lie group acting \textit{linear isometrically} on the real Banach
	space $(X, \|\cdot\|)$ and  $h: X \to R$ be a $G$-invariant, locally Lipschitz functional. 
	If $h|_G$ denotes the restriction of $h$ to ${\sf Fix}_X(G)$ and $u\in {\sf Fix}_X(G)$ is a critical point of $h|_G$ then
	$u$ is also a critical point of $h$.
\end{proposition}
The smooth version of the principle of
symmetric criticality has been provided by Palais \cite{Palais} and later extended to various nonsmooth settings.

\section{Non-existence of solutions: proof of Theorem \ref{theorem-non-existence}}\label{section-3}

In this section we prove  Theorem \ref{theorem-non-existence}.  Let $u\in H^1(M)$ be a solution of \eqref{problem-1}, i.e., relation \eqref{solution-DI} holds for every $v\in H^1(M)$. Let us choose $v=u$ in \eqref{solution-DI}, obtaining 
\begin{equation}\label{v=u}
		\int_M |\nabla_gu(x)|^2{\rm d}v_g -\mu\int_M\frac{u^2(x)}{d^2_g(x_0,x)}{\rm d}v_g+\int_Mu^2(x){\rm d}v_g=\lambda\int_M \alpha(x)\xi_x u(x){\rm d}v_g,
\end{equation} 
where $\xi_x\in \partial F(u(x))$ is a suitable selection,  $x\in M$, such that  $x\mapsto \alpha(x)\xi_x v(x)$ belongs to $L^1(M)$.  By assumptions	$( \textbf{H})_\alpha$,    $(\textbf{H})_0$ and relation \eqref{v=u} we obtain that 
\begin{equation}\label{u=v-modified}
	\int_M |\nabla_gu(x)|^2{\rm d}v_g -\mu\int_M\frac{u^2(x)}{d^2_g(x_0,x)}{\rm d}v_g+\int_Mu^2(x){\rm d}v_g\leq |\lambda|C_0 \|\alpha\|_{L^\infty}\int_Mu^2(x){\rm d}v_g.
\end{equation}
Assume by contradiction that $u\neq 0.$

\underline{Proof of (i)}: ${\bf K}\leq -\kappa$ for some $\kappa\geq 0.$  

Let $\kappa=0$. If $\mu\leq \frac{(n-2)^2}{4}$, by the Hardy inequality \eqref{Hardy-CH}  and relation \eqref{u=v-modified}, it turns out that
$$\int_Mu^2(x){\rm d}v_g< |\lambda|C_0 \|\alpha\|_{L^\infty}\int_Mu^2(x){\rm d}v_g;$$ 
here we used the fact that equality cannot occur in the  Hardy inequality \eqref{Hardy-CH} unless $u=0.$
Consequently, if $|\lambda|C_0\|\alpha\|_{L^\infty}\leq 1,$ we arrive to a contradiction, i.e., we necessarily have $u=0,$ concluding the proof of (i1). 

Let  $\kappa>0$. Assume first that $0< \mu\leq \frac{(n-2)^2}{4}.$ Then by  the Hardy inequality \eqref{Hardy-CH} we have that 
$$\mu\int_M\frac{u^2(x)}{d^2_g(x_0,x)}{\rm d}v_g< \frac{4\mu}{(n-2)^2}\int_M |\nabla_gu(x)|^2{\rm d}v_g,$$
where we again used the fact that no equality occurs in \eqref{Hardy-CH} for non-zero functions. 
Thus, by 
\eqref{u=v-modified} it follows that 
\begin{equation}\label{trivial-1}
	\left(1-\frac{4\mu}{(n-2)^2}\right)\int_M |\nabla_gu(x)|^2{\rm d}v_g< \left(|\lambda|C_0 \|\alpha\|_{L^\infty}-1\right)\int_Mu^2(x){\rm d}v_g.
\end{equation}
First, if $|\lambda|C_0 \|\alpha\|_{L^\infty}\leq 1$, since $ \mu\leq \frac{(n-2)^2}{4},$  relation \eqref{trivial-1} gives a contradiction. Second, if $|\lambda|C_0 \|\alpha\|_{L^\infty}> 1$, by our assumption 
$(n-2)^2(|\lambda|C_0\|\alpha\|_{L^\infty}-1)\leq (n-1)^2\left(\frac{(n-2)^2}{4}-\mu\right)\kappa$ we obtain that $ \mu< \frac{(n-2)^2}{4}$; moreover,  relation \eqref{trivial-1} and the assumption imply that 
$$\int_M |\nabla_gu(x)|^2{\rm
	 d}v_g<\frac{(n-1)^2}{4}\kappa\int_Mu^2(x){\rm d}v_g.$$
The latter inequality is in contradiction to McKean's  spectral gap theorem, see \eqref{McKean-spectral-gap}. Therefore, we necessarily have $u=0$, concluding the proof of (i2) for $\mu> 0$. 

If $\mu\leq 0$, then our assumption reduces to $|\lambda|C_0\|\alpha\|_{L^\infty}-1\leq \frac{(n-1)^2}{4}\kappa$ and 
 by 
 \eqref{u=v-modified} one has that
$$
 \int_M |\nabla_gu(x)|^2{\rm d}v_g\leq \left(|\lambda|C_0 \|\alpha\|_{L^\infty}-1\right)\int_Mu^2(x){\rm d}v_g.
$$ Therefore, we obtain that 
 $$\int_M |\nabla_gu(x)|^2{\rm
 	d}v_g\leq \frac{(n-1)^2}{4}\kappa\int_Mu^2(x){\rm d}v_g.$$
 Since no equality  occurs in  McKean's  spectral gap estimate \eqref{McKean-spectral-gap} for any non-zero function $u\in H^1(M)$, we arrive to a contradiction. In conclusion, we necessarily have that $u=0$, which ends the proof of (i2) also for $\mu\leq 0$.

\underline{Proof of (ii)}: {\rm {\sf Ric}}$_{(M,g)}\geq 0$. Since $\mu\leq {\sf AVR}_{(M,g)}^\frac{2}{n}\frac{(n-2)^2}{4}$,  the Hardy inequality from \eqref{Hardy-Ric} (together with the fact that no non-zero function realizes the equality) and relation \eqref{u=v-modified} imply that 
$$\int_Mu^2(x){\rm d}v_g< |\lambda|C_0 \|\alpha\|_{L^\infty}\int_Mu^2(x){\rm d}v_g.$$ 
Consequently, if $|\lambda|C_0\|\alpha\|_{L^\infty}\leq 1,$ we arrive to a contradiction; thus $u=0.$ This ends the proof of (ii). \hfill $\square$

\section{Sub-quadratic case: proof of Theorem \ref{theorem-multiplicity}}\label{section-4}

Throughout this section we assume the assumptions in Theorem \ref{theorem-multiplicity} are satisfied. The proof is divided into several steps. 

\underline{Step 1.} (Truncation and nonsmooth energy functional) Since $F(0)=0$ (by $({\bf H})_3$), we consider the truncated locally Lipschitz function $F^+(t)=F(t_+)$, $t\in \mathbb R$. The \textit{energy functional} $\mathcal E^+:H^1(M)\to \mathbb R$ to the slightly modified problem \eqref{problem-1}, considering $F^+$ instead of $F$, is defined as 
$$\mathcal E^+(u)=\frac{1}{2}\mathcal N_\mu(u)-\lambda\mathcal F^+(u),$$
where 
$$\mathcal N_\mu(u)=\int_M |\nabla_gu(x)|^2{\rm d}v_g -{\mu}\int_M\frac{u^2(x)}{d^2_g(x_0,x)}{\rm d}v_g+\int_Mu^2(x){\rm d}v_g$$ and $$ \mathcal F^+(u)=\int_M\alpha(x) F^+(u(x)){\rm d}v_g.$$

On one hand, it is clear that $\mathcal N_\mu$ is of class $C^1$ on $H^1(M)$ and due to the Hardy inequalities (i.e., \eqref{Hardy-CH} and \eqref{Hardy-Ric}), for the corresponding values of $\mu$ from the statement of the theorem, $\mathcal N_\mu^{1/2}$ turns out to be equivalent to the usual norm $\|\cdot\|_{H^1}$ on $H^1(M)$, i.e., 
\begin{equation}\label{norm-equivalence}
c_\mu \|u\|_{H^1}^2\leq  \mathcal N_\mu(u)\leq \|u\|_{H^1}^2,\ \ \forall u\in H^1(M),
\end{equation}	where 
$$0<c_\mu=
\begin{cases} 
	 1-\frac{4\mu}{(n-2)^2}, & \text{in the case (i)};  \\
	1-{\sf AVR}_{(M,g)}^{-\frac{2}{n}}\frac{4\mu}{(n-2)^2} & \text{in the case (ii)}. 
\end{cases}
$$

 On the other hand, one can prove that $\mathcal F^+$ is well-defined and locally Lipschitz on $H^1(M)$. To see this, we first observe that by  $(\textbf{H})_1$ and $(\textbf{H})_2$, for every $\epsilon>0$ there exists $\delta_\epsilon\in (0,1)$ such that
  \begin{equation}\label{eps-egyenlotlenseg}
  	|\xi|\leq \epsilon t,\ \forall \xi\in \partial F^+(t),\ \forall 0<t<\delta_\epsilon \ \&\ t>\delta_\epsilon^{-1}. 
  \end{equation}
Fix $\epsilon_0>0.$ Since $\partial F^+$ is an upper semicontinuous set-valued map with non-empty compact values, see Proposition \ref{prop-lok-lip-0}/$(i)$, we also have for some $K_{\epsilon_0}>0$  that $|\xi|\leq K_{\epsilon_0} t$ for every $\xi\in   \partial F^+(t)$ and $t\in [\delta_{\epsilon_0},\delta_{\epsilon_0}^{-1}]$. The latter fact with  \eqref{eps-egyenlotlenseg} implies that $$|\xi|\leq C_{\epsilon_0} t,\ \forall \xi\in \partial F^+(t),\ \forall t>0,$$
 where $C_{\epsilon_0}=\max\{{\epsilon_0},K_{\epsilon_0}\}.$ Now, let $u\in H^1(M)$ and $U_u$ be any open bounded neighborhood of $u$ in $H^1(M)$, i.e., for some $K>0$  we have $\|w\|_{H^1}\leq K$ for every $w\in U_u$. If $u_1,u_2\in U_u$, then by Lebourg's mean value theorem, see Proposition \ref{prop-lok-lip-0}/$(iv)$, for a.e.\ $x\in M$ there exist $\gamma\in [0,1]$ and $\xi_x^\gamma\in \partial F^+((1-\gamma)u_1(x)+\gamma u_2(x))$ such that 
  $$|F^+(u_1(x))-F^+(u_2(x))|=|\xi_x^\gamma||u_1(x)-u_2(x)|\leq C_{\epsilon_0}(|u_1(x)|+|u_2(x)|)|u_1(x)-u_2(x)|.$$
  By H\"older's inequality and the trivial embedding $H^1(M)\subset L^2(M)$, we have that 
  $$|\mathcal F^+(u_1)-\mathcal F^+(u_2)|\leq \int_M\alpha(x) |F^+(u_1(x))-F^+(u_2(x))|{\rm d}v_g\leq 2C_{\epsilon_0}\|\alpha\|_{L^\infty}K\|u_1-u_2\|_{H^1},$$
  which means that $\mathcal F^+$ is Lipschitz on $U_u$. The fact that $\mathcal F^+$ is well-defined follows in a similar way.  Having these properties, a similar argument as in   Clarke \cite[Section 2.7]{Clarke} (see also Costea, Krist\'aly and Varga \cite{CKV}) shows that for every closed subspace $W$ of $H^1(M)$ we have that  $$\partial (\mathcal F^+|_W)(u)\subseteq \int_M \alpha(x)\partial F^+(u(x)){\rm d}v_g,\ \ \forall u\in W;$$
here, $\mathcal F^+|_W$ is the restriction of the functional $\mathcal F^+$ to the subspace $W$ and the latter inclusion has the following interpretation: to every $\xi\in \partial (\mathcal F^+|_W)(u)$ there exists a measurable selection $x\mapsto \xi_x\in \partial F^+(u(x))$ such that  the map $x\mapsto \alpha(x)\xi_x w(x)$ belongs to $L^1(M)$ for every $w\in W$ and 
$$\langle\xi,w \rangle=\int_M \alpha(x)\xi_xw(x){\rm d}v_g.$$
By using Fatou's lemma, Lebourg's mean value theorem, Lebesgue's  dominated convergence theorem, and a careful limiting argument, see e.g. Krist\'aly \cite{Kristaly-SVAN} in the Euclidean setting, it turns out that 
\begin{equation}\label{Rocky-type}
	 (\mathcal F^+|_W)^0(u;w)\leq \int_M\alpha(x) (F^+)^0(u(x);w(x)){\rm d}v_g,\ \ \forall u,w\in W.
\end{equation}

Let $u\in H^1(M)$ be a critical point of $\mathcal E^+$, i.e., $0\in \partial \mathcal E^+(u)$. We are going to prove that $u$ is a non-negative solution to the differential inclusion \eqref{problem-1}. First, by Proposition \ref{prop-lok-lip-0}/$(v)\&(vii)$, we have that $$\frac{1}{2}\mathcal N_\mu'(u)\in \lambda\partial \mathcal F^+(u),$$ i.e., for every test-function $w\in H^1(M)$ one has
$$	\int_M \nabla_gu(x)\nabla_gw(x){\rm d}v_g -\mu\int_M\frac{u(x)w(x)}{d^2_g(x_0,x)}{\rm d}v_g+\int_Mu(x)w(x){\rm d}v_g=\lambda\int_M \alpha(x)\xi_x w(x){\rm d}v_g,$$
with the above interpretation for the right hand side. Let $u_-=\min(0,u)$ be the non-positive part of $u$ and note that it belongs to the space $H^1(M)$, see Hebey \cite[Proposition 2.5]{Hebey}. If we put $v=u_-$ into the latter relation, we obtain that $\xi_x u_-(x)=0$ for a.e. $x\in M$ since $\xi_x\in \partial F^+(u(x))$ (thus $\xi_x=0 $ whenever $u(x)<0$). In consequence, $\mathcal N_\mu(u_-)=0$, thus $u_-=0$, i.e., $u\geq 0$. In particular,  $\xi_x\in \partial F^+(u(x))=\partial F(u(x))$, therefore the latter relation is precisely \eqref{solution-DI},  which means that $u\in H^1(M)$ is a non-negative solution of  \eqref{problem-1}. 

\underline{Step 2.} (Isometry actions) Let $G$ be a compact connected subgroup of ${\sf Isom}_g(M )$ with the property that ${\sf Fix}_M(G) =\{x_0\}$ for the same $x_0\in M$ as in problem \eqref{problem-1}. The action of $G$ on $H^1(M)$, i.e., $G\times H^1(M)\to H^1(M)$, is defined by 
\begin{equation}\label{action}
	(\sigma u)(x)=u(\sigma^{-1}(x)),\ \forall \sigma\in G,\ u\in H^1(M),\ x\in M.
\end{equation}
It is standard to prove that  $G$ acts continuously and linearly on $H^1(M)$. For instance, if $\sigma_1,\sigma_2\in G$, it turns out that for every $u\in H^1(M)$ and $\sigma\in G$ we have
$$(\sigma_1\circ \sigma_2)u(x)=u((\sigma_1\circ \sigma_2)^{-1}(x))=u(\sigma_2^{-1}( \sigma_1^{-1}(x)))=(\sigma_2u)(\sigma_1^{-1}(x))=(\sigma_1(\sigma_2u))(x).$$ Moreover, since $G$ contains  isometries of $(M,g)$, the functionals  $\displaystyle u\mapsto \int_M |\nabla_gu(x)|^2{\rm d}v_g$ and  $\displaystyle u\mapsto \int_M u^2(x){\rm d}v_g$ are both $G$-invariant; in particular,  $\|\sigma u\|_{H^1}=\| u\|_{H^1}$ for every $\sigma\in G$ and $u\in H^1(M)$, i.e., $G$ acts isometrically on $H^1(M).$

Since  ${\sf Fix}_M(G) =\{x_0\}$, it turns out that for every $\sigma\in G$ and $y\in M$, we have $d_g(x_0,\sigma(y))=d_g(\sigma(x_0),\sigma(y))=d_g(x_0,y)$; therefore, 
 a change of variables implies that
\begin{eqnarray*}
	\int_M\frac{(\sigma u)^2(x)}{d^2_g(x_0,x)}{\rm d}v_g(x)&=&\int_M\frac{u^2(\sigma^{-1}(x))}{d^2_g(x_0,x)}{\rm d}v_g(x)=\int_M\frac{u^2(y)}{d^2_g(x_0,\sigma(y))}{\rm d}v_g(\sigma(y))\\&=&\int_M\frac{u^2(y)}{d^2_g(x_0,y)}{\rm d}v_g(y).
\end{eqnarray*}
In particular, the functional $u\mapsto \mathcal N_\mu(u)$ is $G$-invariant on $H^1(M)$.  

Furthermore, since $\alpha:M\to \mathbb R$  depends only on $d_g(x_0,\cdot)$, it  is also $G$-invariant, and one can prove by a change of variables that for every $\sigma\in G$ and $u\in H^1(M)$, 
 \begin{eqnarray*}
 \mathcal F^+(\sigma u)&=&\int_M\alpha(x) F^+((\sigma u)(x)){\rm d}v_g(x)=\int_M\alpha(x) F^+(u(\sigma^{-1}(x))){\rm d}v_g(x)\\&=&\int_M\alpha(\sigma(y)) F^+(u(y)){\rm d}v_g(\sigma(y))=\int_M\alpha(y) F^+(u(y)){\rm d}v_g(y)\\&=&
 \mathcal F^+(u), \end{eqnarray*} 
i.e., $\mathcal F^+$ is $G$-invariant on $H^1(M)$. 
In conclusion, the energy functional $\mathcal E^+=\mathcal N_\mu/2-\lambda\mathcal F^+$ is $G$-invariant on $H^1(M)$. 

The set of fixed points  of $G$ over $H^1(M)$, i.e., ${\sf Fix}_{H^1(M)}(G)$, is nothing but the closed subset of $G$-invariant functions of $H^1(M)$. Now, according to the principle of symmetric criticality, see Proposition \ref{psc}, 
if $u\in {\sf Fix}_{H^1(M)}(G)=:H^1_G(M)$ is a critical point of the restricted energy functional $\mathcal E^+_G:=\mathcal E^+|_{H_G^1(M)}$ then $u$ is also a critical point of the initial energy functional $\mathcal E^+$.  

\underline{Step 3.} (Spectral gap estimate for $\mathcal F^+/\mathcal N_\mu$ on $H_G^1(M)$.) 
We are going to prove that for every admissible $\mu$ from the statement of the theorem, one has
\begin{equation}\label{fontos-suprem}
	0<\sup_{u\in H_G^1(M)\setminus \{0\}}\frac{\mathcal F^+(u)}{\mathcal N_\mu(u)}<+\infty.
\end{equation}

Let $q\in (2,2^*)$ and fix arbitrarily $\epsilon>0$ together with the number $\delta_{\epsilon}>0$  appearing in \eqref{eps-egyenlotlenseg}. By the boundedness of the function $t\mapsto \frac{\max|\partial F^+(t)|}{t^{q-1}}$ on $[\delta_{\epsilon},\delta_{\epsilon}^{-1}]$ and due to \eqref{eps-egyenlotlenseg}, there exists $l_\epsilon>0$ such that 
\begin{equation}\label{eps-elso}
	0\leq |\xi|\leq \epsilon t+l_\epsilon t^{q-1},\ \forall t\geq 0,\ \xi\in \partial F^+(t)=\partial F(t).
\end{equation}
Note that we have
\begin{equation}\label{kell-hatulra}
	0\leq |F^+(t)|\leq \epsilon t^2+l_\epsilon |t|^{q},\ \forall t\in \mathbb R.
\end{equation}
Indeed, on one hand, by definition $F^+(t)=F(0)=0$ for every $t\leq 0$, thus the latter relation trivially holds. On the other hand, for $t\geq 0$, the estimate \eqref{eps-elso}  and Lebourg's mean value theorem immediately imply the required estimate. 

Consequently, the  estimate \eqref{kell-hatulra} shows that for every $u\in H_G^1(M)$ we have
\begin{eqnarray*}
	0\leq |\mathcal F^+(u)|&=&\left|\int_M\alpha(x) F^+(u(x)){\rm d}v_g\right|\leq \int_M\alpha(x) \left|F^+(u(x))\right|{\rm d}v_g\\&\leq & \|\alpha\|_{L^\infty}\left(\epsilon \|u\|_{H^1}^2+l_\epsilon (K_q^\pm)^q\|u\|_{H^1}^{q}\right),
\end{eqnarray*}
where $K_q^\pm>0$ are the embedding constants from \eqref{Sobolev-CH} and 
\eqref{Sobolev-Ric}, respectively. Accordingly, for every $u\in H_G^1(M)\setminus \{0\}$ one has  that
$$0\leq \frac{|\mathcal F^+(u)|}{\mathcal N_\mu(u)}\leq 
c_\mu^{-1}\|\alpha\|_{L^\infty}\left(\epsilon +l_\epsilon (K_q^\pm)^q\|u\|_{H^1}^{q-2}\right),$$
where $c_\mu>0$ is the constant from \eqref{norm-equivalence}. Due to the fact that $q>2$ and $\epsilon>0$ is arbitrarily fixed, it turns out that 
\begin{equation}\label{limit-nulla}
\frac{\mathcal F^+(u)}{\mathcal N_\mu(u)}\to 0\ \ {\rm as}\ \ 	\|u\|_{H^1}\to 0, u\in H_G^1(M).
\end{equation}

The counterpart of \eqref{limit-nulla} at 'infinity' reads as
 \begin{equation}\label{limit-vegtelen}
 	\frac{\mathcal F^+(u)}{\mathcal N_\mu(u)}\to 0\ \ {\rm as}\ \ 	\|u\|_{H^1}\to +\infty, u\in H_G^1(M).
 \end{equation}
Indeed, combining  the boundedness of $t\mapsto \frac{\max|\partial F^+(t)|}{t^{1/2}}$  on $[\delta_{\epsilon},\delta_{\epsilon}^{-1}]$ with the estimate \eqref{eps-egyenlotlenseg}, one can find  $L_\epsilon>0$ such that 
\begin{equation}\label{eps-masodik}
	0\leq |\xi|\leq \epsilon t+L_\epsilon t^{1/2},\ \forall t\geq 0,\ \xi\in \partial F^+(t)=\partial F(t).
\end{equation}
Due to hypothesis $(\textbf{H})_\alpha,$ one has that $\alpha\in L^4(M)$. Then using Lebourg's mean value theorem and H\"older's inequality, we can proceed as before, obtaining 
\begin{equation}\label{becsles-vegtelen}
		0\leq |\mathcal F^+(u)|\leq \int_M\alpha(x) \left|F^+(u(x))\right|{\rm d}v_g\leq  \epsilon\|\alpha\|_{L^\infty} \|u\|_{H^1}^2+L_\epsilon \|\alpha\|_{L^4}\|u\|_{H^1}^\frac{3}{2}.
\end{equation}
Consequently, for every $u\in H_G^1(M)\setminus \{0\}$ we have
$$0\leq \frac{|\mathcal F^+(u)|}{\mathcal N_\mu(u)}\leq 
c_\mu^{-1}\left(\epsilon\|\alpha\|_{L^\infty} +L_\epsilon \|\alpha\|_{L^4}\|u\|_{H^1}^{-\frac{1}{2}}\right).$$
This estimate together with the arbitrariness of $\epsilon>0$ immediately imply \eqref{limit-vegtelen}. 

In particular, \eqref{limit-nulla} and \eqref{limit-vegtelen} imply that the second inequality in \eqref{fontos-suprem} holds. In order to check the first inequality in \eqref{fontos-suprem}, we recall by $(\textbf{H})_3$ that there exists $t_0^+>0$ such that $F(t_0^+)>0.$ Moreover, by  $(\textbf{H})_\alpha,$ since $\alpha\neq 0$ and it depends only on $d_g(x_0,\cdot)$, there exists an open $x_0$-centered annulus on $M$ with radii $0\leq r<R$, i.e. $A_{x_0}(r,R)=\{x\in M:r<d_g(x_0,x)<R\}$, such that
${\rm essinf}_{A_{x_0}(r,R)}\alpha=\alpha_0>0$. For sufficiently small $\epsilon>0$ (e.g. $\epsilon<(R-r)/3$), we  consider the function $w_\epsilon:M\to \mathbb R$ defined by  $$w_\epsilon(x)=\left\{ \begin{array}{lll}
	\frac{t_0^+}{\epsilon}(d_g(x_0,x)-r) &\mbox{if} &  d_g(x_0,x)\in (r,r+\epsilon),\\
	t_0^+ &\mbox{if} &  d_g(x_0,x)\in [r+\epsilon,R-\epsilon],\\
		\frac{t_0^+}{\epsilon}(R-d_g(x_0,x)) &\mbox{if} &  d_g(x_0,x)\in (R-\epsilon,R),\\
			0 &\mbox{if} & x\notin A_{x_0}(r,R).
\end{array}\right.$$
Note that $w_\epsilon\in H_G^1(M)$ and $w_\epsilon\geq 0$. Moreover, 
\begin{eqnarray*}
	\mathcal F^+(w_\epsilon)&=&\int_M\alpha(x) F(w_\epsilon(x)){\rm d}v_g=\int_{A_{x_0}(r,R)}\alpha(x) F(w_\epsilon(x)){\rm d}v_g\\&\geq &
	\alpha_0 F(t_0^+){V}_g(A_{x_0}(r+\epsilon,R-\epsilon))\\&&-\|\alpha\|_{L^\infty}\max_{t\in [0,t_0^+]}|F(t)|(V_g(A_{x_0}(r,r+\epsilon))+V_g(A_{x_0}(R-\epsilon,R))).
\end{eqnarray*}
By continuity reason,  there exists $\epsilon_0>0$ such that for every $\epsilon\in (0,\epsilon_0)$, 
$$\mathcal F^+(w_\epsilon)\geq \alpha_0 F(t_0^+){V}_g(A_{x_0}(r,R))/2>0.$$
On the other hand, by \eqref{norm-equivalence} and the eikonal equation ($|\nabla_g d_g(x_0,\cdot)|=1$ a.e.\ on $M$) we have the  estimate 
$$\mathcal N_\mu(w_\epsilon)\leq \|w_\epsilon\|_{H^1}^2\leq (t_0^+)^2(1+\epsilon^{-2}){V}_g(A_{x_0}(r,R))<+\infty.$$
Consequently, it turns out that
$$0<\frac{\mathcal F^+(w_{\epsilon_0/2})}{\mathcal N_\mu(w_{\epsilon_0/2})}\leq \sup_{u\in H_G^1(M)\setminus \{0\}}\frac{\mathcal F^+(u)}{\mathcal N_\mu(u)},$$
which shows the validity of the first inequality in \eqref{fontos-suprem}. 

\underline{Step 4.} (Analytic properties of $\mathcal E_G^+$) We shall prove three basic properties of $\mathcal E_G^+$ on $H^1_G(M)$, namely,  coercivity and boundedness from below, as well as the validity of the nonsmooth Palais-Smale condition. 

Let $\lambda>0$ be arbitrarily fixed and $\mu$ be in the admissible range (cf.\ the statement of the theorem). First, we observe by \eqref{norm-equivalence} and \eqref{becsles-vegtelen} that 
for every $u\in H^1_G(M)$ we have
\begin{eqnarray*}
	\mathcal E_G^+(u)&=&\frac{1}{2}\mathcal N_\mu(u)-\lambda\mathcal F^+(u)\\
	&\geq&\left(\frac{c_\mu}{2}-\epsilon\lambda\|\alpha\|_{L^\infty}\right)   \|u\|_{H^1}^2-\lambda L_\epsilon \|\alpha\|_{L^4}\|u\|_{H^1}^\frac{3}{2}. 
\end{eqnarray*}
In particular, for sufficiently small $\epsilon>0$, e.g.\ 
$0<\epsilon<\frac{c_\mu}{2}\lambda^{-1}\|\alpha\|_{L^\infty}^{-1}$, it follows that $\mathcal E_G^+$ is bounded from below and coercive, i.e., $\mathcal E_G^+(u)\to +\infty$ whenever $\|u\|_{H^1}\to +\infty.$

Let $\{u_k\}_k\subset H^1_G(M)$ be a Palais-Smale sequence for $\mathcal E_G^+$, i.e., for some $M>0$, one has  $|\mathcal E_G^+(u_k)|\leq M$ and $m(u_k)\to 0$ as $k\to \infty$, where $m(u)=\min\{\|\xi\|_*:\xi\in \partial \mathcal E_G^+(u)\}$. We want to prove that, up to a subsequence, $\{u_k\}_k$ strongly converges to some element in $H^1_G(M)$. 
Being $\mathcal E_G^+$ coercive, the sequence  $\{u_k\}_k\subset H^1_G(M)$ is bounded in $H^1_G(M)$. Therefore, due to the fact that $H_G^1(M)$ can be compactly embedded into $L^q(M)$, $q\in (2,2^*)$, see  \S \ref{subsection-compact},  up to a subsequence, one has that 
\begin{equation}\label{u-k-weakly}
	u_k\to u\ {\rm weakly\ in}\ H_G^1(M);
\end{equation}
\begin{equation}\label{u-k-strongly}
	u_k\to u\ {\rm strongly\ in}\ L^q(M), \ q\in (2,2^*).
\end{equation}
By Proposition \ref{prop-lok-lip-0}/$(v)$ and the definition of $\mathcal E_G^+$ we have that
$$(\mathcal E_G^+)^0(u_k;u-u_k)=\frac{1}{2}\langle \mathcal N_\mu'(u_k),u-u_k\rangle + \lambda(-\mathcal F^+)^0(u_k;u-u_k);$$
$$(\mathcal E_G^+)^0(u;u_k-u)=\frac{1}{2}\langle \mathcal N_\mu'(u),u_k-u\rangle +\lambda(-\mathcal F^+)^0(u;u_k-u).$$
Note that 
$$\frac{1}{2}\langle \mathcal N_\mu'(u_k),u-u_k\rangle+\frac{1}{2}\langle \mathcal N_\mu'(u),u_k-u\rangle=-\mathcal N_\mu(u_k-u).$$
By  adding the above relations and using Proposition \ref{prop-lok-lip-0}/$(vi)$, it turns out that
\begin{eqnarray}\label{summative-relation}
\nonumber	\mathcal N_\mu(u_k-u)&=& \lambda \left((\mathcal F^+)^0(u_k;-u+u_k)+(\mathcal F^+)^0(u;-u_k+u)\right)\\&&-(\mathcal E_G^+)^0(u_k;u-u_k)-(\mathcal E_G^+)^0(u;u_k-u).
\end{eqnarray}
In the sequel, we are going to estimate the terms in the right hand side of \eqref{summative-relation}. First, by inequality \eqref{Rocky-type},  Proposition \ref{prop-lok-lip-0}/$(ii)$ and \eqref{eps-elso} together with the fact that $\partial F^+(t)=\{0\}$ for $t\leq 0$, we have
\begin{eqnarray*}
	I_k^1&:=&(\mathcal F^+)^0(u_k;-u+u_k)+(\mathcal F^+)^0(u;-u_k+u)\\&\leq&
	\int_M\alpha(x)\left[ (F^+)^0(u_k(x);u_k(x)-u(x))+(F^+)^0(u(x);u(x)-u_k(x))\right]{\rm d}v_g\\&= &	\int_M\alpha(x)\left[ \max \{\xi_k(u_k(x)-u(x)):\xi_k\in \partial F^+(u_k(x))\}\right.\\&&\left.\qquad\qquad + \max \{\xi(u(x)-u_k(x)):\xi\in \partial F^+(u(x))\}\right]{\rm d}v_g
		\\&\leq &\|\alpha\|_{L^\infty}\int_M[\epsilon(|u_k(x)|+|u(x)|)+l_\epsilon(|u_k(x)|^{q-1}+|u(x)|^{q-1})]|u(x)-u_k(x)|{\rm d}v_g\\&\leq&
	2\epsilon\|\alpha\|_{L^\infty}(\|u_k\|^2_{H^1}+\|u\|^2_{H^1})+l_\epsilon\|\alpha\|_{L^\infty}(\|u_k\|_{L^q}^{q-1}+\|u\|^{q-1}_{L^q})\|u_k-u\|_{L^q}.
\end{eqnarray*}
By the arbitrariness of $\epsilon>0$ and the convergence property \eqref{u-k-strongly}, the latter estimate shows that 
\begin{equation}\label{I_1}
	\limsup_{k\to \infty}	I_k^1\leq 0. 
\end{equation}

Let $\xi_k\in \partial \mathcal E_G^+(u_k)$ be such that $m(u_k)=\|\xi_k\|_*$. Thus, we have that
\begin{eqnarray*}
	I_k^2:=(\mathcal E_G^+)^0(u_k;u-u_k)\geq \langle \xi_k,u-u_k\rangle\geq -\|\xi_k\|_*\|u-u_k\|_{H^1}. 
\end{eqnarray*}
Consequently, since $m(u_k)=\|\xi_k\|_*\to 0$ as $k\to \infty$, we have that
\begin{equation}\label{I_2}
	\liminf_{k\to \infty}	I^2_k\geq 0. 
\end{equation}

Moreover, for every $\xi\in \partial \mathcal E_G^+(u)$, we also have that
$I^3_k:=(\mathcal E_G^+)^0(u;u_k-u)\geq \langle \xi,u_k-u\rangle$; thus, by the weak convergence property \eqref{u-k-weakly} we have that 
\begin{equation}\label{I_3}
	\liminf_{k\to \infty}	I^3_k\geq 0. 
\end{equation}

By the estimates \eqref{I_1}-\eqref{I_3} and relation \eqref{summative-relation} we have that 
$$0\leq \limsup_{k\to \infty}	\mathcal N_\mu(u_k-u)\leq  \limsup_{k\to \infty}I^1_k-\liminf_{k\to \infty}I^2_k-\liminf_{k\to \infty}I^3_k\leq 0,$$
i.e., 	$\mathcal N_\mu(u_k-u)\to 0$ as $k\to \infty.$ Due to \eqref{norm-equivalence}, it turns out that $u_k\to u$ strongly in the $H^1$-norm as $k\to \infty$, which is the desired property.

\underline{Step 5.} (Local minimum point for $\mathcal E_G^+$: first solution) 
Let $$\lambda_0^+:=\inf_{\substack{u\in H_G^1(M) \\
		\mathcal F^+(u)>0}}\frac{\mathcal N_\mu(u)}{2\mathcal F^+(u)}.$$
Due to Step 3, see \eqref{fontos-suprem},  one has that $0<\lambda_0^+<\infty.$

 If we fix $\lambda> \lambda_0^+$, one can find $\tilde w_\lambda\in H_G^1(M)$ with $\mathcal F^+(\tilde
w_\lambda)>0$ such that $$\lambda>\frac{\mathcal N_\mu(\tilde
	w_\lambda)}{2\mathcal F^+(\tilde
	w_\lambda)}\geq \lambda_0^+.$$
Thus, by the latter inequality we have
$$C_\lambda^1:=\inf_{H_G^1(M)}\mathcal E_G^+\leq \mathcal E_G^+(\tilde
w_\lambda)=\frac{1}{2}\mathcal N_\mu(\tilde
w_\lambda)-\lambda\mathcal F^+(\tilde
w_\lambda)< 0.$$
Due to the fact that  $\mathcal E_G^+$ is bounded from below and verifies the nonsmooth
Palais-Smale condition (see Step 4),  $C_\lambda^1$ is a critical value
of $\mathcal E_G^+$, see Chang \cite[Theorem 3.5]{Chang},  i.e.,  there exists $u_\lambda^1\in
H^1_G(M)$ such that $\mathcal
E_G^+(u_\lambda^1)=C_\lambda^1<0$ and $0\in \partial \mathcal
E_G^+(u_\lambda^1).$ In particular, $u_\lambda^1\neq 0$ (since $\mathcal
E_G^+(u_\lambda^1)<0=\mathcal
E_G^+(0)$), and by the principle of symmetric criticality, $u_\lambda^1$ is a critical point also for the initial energy functional (see Step 2), i.e., $0\in \partial \mathcal
E^+(u_\lambda^1).$ According to (the final part of) Step 1, $u_\lambda^1\in H_G^1(M)$ is a non-negative solution to the differential inclusion \eqref{problem-1}. 

\underline{Step 6.} (Minimax-type critical point for $\mathcal E_G^+$: second solution) Let $\lambda> \lambda_0^+$.
Due to  \eqref{kell-hatulra}, for  sufficiently small $\epsilon>0$   (e.g., $\frac{c_\mu}{4}>\epsilon\lambda\|\alpha\|_{L^\infty}$)
and for every $u\in H_G^1(M)$ one has that 
\begin{eqnarray*}
	\mathcal E_G^+(u) &=& \frac{1}{2}\mathcal N_\mu(u)-\lambda\mathcal F^+(u)\geq \left(\frac{c_\mu}{2}-\epsilon\lambda\|\alpha\|_{L^\infty}\right)   \|u\|_{H^1}^2-\lambda \|\alpha\|_{L^\infty} l_\epsilon (K_q^\pm)^q\|u\|_{H^1}^{q},
\end{eqnarray*}
where $q\in (2,2^*)$ and $K_q^\pm>0$ are the embedding constants from \eqref{Sobolev-CH} and 
\eqref{Sobolev-Ric}, respectively. 
Let
$$\rho_\lambda=\min\left\{\|\tilde w_\lambda\|_{H^1},\left(
\frac{\frac{c_\mu}{2}-\epsilon\lambda\|\alpha\|_{L^\infty}}{2\lambda \|\alpha\|_{L^\infty} l_\epsilon (K_q^\pm)^q}\right)^\frac{1}{q-2}
\right\}.$$ The choice of $\rho_\lambda>0$ and Step 4 show that
$$\inf_{\|u\|_{H^1}=\rho_\lambda; u\in H_G^1(M)}\mathcal E_G^+(u)\geq \frac{1}{2}\left(\frac{c_\mu}{2}-\epsilon\lambda\|\alpha\|_{L^\infty}\right)\rho_\lambda^2>0=\mathcal E_G^+(0)> \mathcal E_G^+(\tilde
w_\lambda).$$ The latter estimate shows that the functional $\mathcal E_G^+$ has the
 mountain pass geometry. On account of Step 4, since $\mathcal E_G^+$ satisfies the nonsmooth Palais-Smale condition, we can  apply
the mountain pass theorem for locally Lipschitz functions, see e.g.\ Kourogenis and Papageorgiou \cite{KP} or Krist\'aly,  Motreanu and  Varga \cite[Theorem 2]{KMV}, guaranteeing the existence of  $u_\lambda^2\in
H_G^1(M)$ with the properties  $0\in \partial \mathcal E_G^+(u_\lambda^2)$ and
$$\mathcal E_G^+(u_\lambda^2)=C_\lambda^2=\inf_{\gamma\in \Gamma}\max_{t\in [0,1]}\mathcal E_G^+(\gamma(t)),$$
where $$\Gamma=\{\gamma\in C([0,1];H_G^1(M)):\gamma(0)=0,\
\gamma(1)=\tilde w_\lambda\}.$$ Since $$C_\lambda^2\geq
\inf_{\|u\|_{H^1}=\rho_\lambda; u\in H_G^1(M)}\mathcal E_G^+(u)>0,$$ it is
clear that $0\neq u_\lambda^2\neq u_\lambda^1.$ The rest of the proof is similar to the end of Step 5, which shows that 
$u_\lambda^2\in H_G^1(M)$ is a non-negative solution to the differential inclusion \eqref{problem-1}, different from $u_\lambda^1$. 

\underline{Step 7.} (Repetition of Steps 1-6 for $\mathcal E^-$) Let $F^-(t)=F(t_-)$, $t\in \mathbb R$, where $t_-=\min(t,0)$. The  locally Lipschitz energy functional $\mathcal E^-:H^1(M)\to \mathbb R$ is defined as
$$\mathcal E^-(u)=\frac{1}{2}\mathcal N_\mu(u)-\lambda\mathcal F^-(u),$$
where 
 $$ \mathcal F^-(u)=\int_M\alpha(x) F^-(u(x)){\rm d}v_g.$$
One can show that if $u\in H^1(M)$ is a critical point of $\mathcal E^-$, i.e., $0\in \partial \mathcal E^-(u)$, then 
it is a non-positive solution of  \eqref{problem-1}, cf. Step 1. 

By using the isometry action \eqref{action}, one can prove that $\mathcal F^-$ is $G$-invariant on $H_G^1(M)$, and if $u\in {\sf Fix}_{H^1(M)}(G)=:H^1_G(M)$ is a critical point of  $\mathcal E^-_G:=\mathcal E^-|_{H_G^1(M)}$ then $0\in \partial \mathcal E^-(u)$ as well, cf.\ Step 2. 

 Instead of the spectral gap estimate \eqref{fontos-suprem}, one can prove 
$$
	0<\sup_{u\in H_G^1(M)\setminus \{0\}}\frac{\mathcal F^-(u)}{\mathcal N_\mu(u)}<+\infty,
$$
cf.\ Step 3, and similar analytic properties are valid for $\mathcal E_G^-$ as in Step 4 (i.e, coercivity,  boundedness from below, and the validity of the nonsmooth Palais-Smale condition). Here, we use again the  compact embedding results from   \S \ref{subsection-compact}. 

Finally, if
$$\lambda_0^-:=\inf_{\substack{u\in H_G^1(M) \\
		\mathcal F^-(u)>0}}\frac{\mathcal N_\mu(u)}{2\mathcal F^-(u)},$$
by the previous part we know that $0<\lambda_0^-<\infty$ and similarly to Steps 5 and 6, we can guarantee for every $\lambda>\lambda_0^-$ a local minimum point $u_\lambda^3\in H_G^1(M)$ of $\mathcal E_G^-$ with $\mathcal E_G^-(u_\lambda^3)<0$ and a minimax-type point $u_\lambda^4\in H_G^1(M)$ of $\mathcal E_G^-$ with $\mathcal E_G^-(u_\lambda^4)>0$; in particular, $u_\lambda^3\neq u_\lambda^4$ and none of them is trivial. These elements are $G$-invariant, non-positive solutions to the differential inclusion \eqref{problem-1}. 

If we choose $\lambda_0=\max(\lambda_0^+,\lambda_0^-)$, we can apply the above arguments, providing four different, non-zero  $G$-invariant solutions  to the differential inclusion \eqref{problem-1} for every $\lambda>\lambda_0$, two of them being non-negative and the other two being non-positive. The proof is complete. 
\hfill $\square$

\section{Super-quadratic case: proof of Theorem \ref{theorem-multiplicity-2}}\label{section-5}

We assume in the sequel that the assumptions of  Theorem \ref{theorem-multiplicity-2} are fulfilled. We again divide the proof into some steps.

\underline{Step 1.} (Functional setting) In view of the previous section, this part is standard. Indeed, the energy functional $\mathcal E:H^1(M)\to \mathbb R$ is defined as
$$\mathcal E(u)=\frac{1}{2}\mathcal N_\mu(u)-\lambda\mathcal F(u),$$
where 
$$ \mathcal F(u)=\int_M\alpha(x) F(u(x)){\rm d}v_g.$$

 Note that by  $(\textbf{H})_1$ and $(\textbf{H})_5$, for every $\epsilon>0$ there exists $C_\epsilon>0$ such that
 \begin{equation}\label{eps-harmadik}
  |\xi|\leq \epsilon |t|+C_\epsilon |t|^{q-1},\ \forall t\in \mathbb R,\ \xi\in \partial F(t).
 \end{equation}
Consequently, 
one has
 \begin{equation}\label{eps-harmadik-bis}
|F(t)|\leq \epsilon t^2+C_\epsilon |t|^{q},\ \forall t\in \mathbb R.
\end{equation}
Since $2<q<2+\frac{4}{n}<2^*$, by using Lebourg's mean value theorem and \eqref{eps-harmadik}, one can prove that $\mathcal F$ is well-defined and locally Lipschitz on $H^1(M)$. It is now standard to show that any critical point $u\in H^1(M)$  of $\mathcal E$ is  a solution of  \eqref{problem-1}.

\underline{Step 2.} (Isometry actions) By using the action \eqref{action}, one can prove in a similar way as in \S \ref{section-4} that $\mathcal E$ is $G$-invariant on $H^1(M).$ Moreover,  the principle of symmetric criticality (Proposition \ref{psc}) implies that if $u\in {\sf Fix}_{H^1(M)}(G)=:H^1_G(M)$ is a critical point of  $\mathcal E_G:=\mathcal E|_{H_G^1(M)}$ then $u$ is also a critical point of $\mathcal E$. 
 
\underline{Step 3.} (Super-quadracity of $F$ at infinity) We are going to prove that 
\begin{equation}\label{super-quadratic}
	F(t)\geq \frac{C}{\nu-2}|t|^\nu,\ \forall t\in \mathbb R,
\end{equation} 
where $\nu>2$ and $C>0$ come from hypothesis $(\textbf{H})_4;$ this means in particular that $F$ is  super-quadratic at infinity (as $\nu>2$). To do this, let 
$h:\mathbb R\to \mathbb R$ be defined by 
$$h(t)=t^{-2}F(t)-\frac{C}{\nu-2}|t|^{\nu-2},\ t\neq 0,$$
and $h(0)=0$. Note that $h$ is well-defined  and locally Lipschitz (indeed, by $(\textbf{H})_1$ and $F(0)=0$ we have that $F(t)=o(t^2)$ as $t\to 0$). By Proposition \ref{prop-lok-lip-0}/$(v)$, one has 
$$\partial h(t)=-2t^{-3}F(t)+t^{-2}\partial F(t)-C|t|^{\nu-4}t,\ \forall t\in \mathbb R\setminus \{0\}.$$

We shall prove \eqref{super-quadratic} for $t\geq 0$, the case $t\leq 0$ being similar. Let $t> 0$; then by Lebourg's mean value theorem, there exist $\theta\in (0,t)$ and $\xi_h\in \partial h(\theta)$ such that
$h(t)=h(t)-h(0)=\xi_ht$. In turn, there exists $\xi_F\in \partial F(\theta)$ such that
$\xi_h=-2\theta^{-3}F(\theta)+\theta^{-2}\xi_F-C\theta^{\nu-3}$ and by  $(\textbf{H})_4$ we have that
\begin{eqnarray*}
h(t)&=&\xi_ht=(-2\theta^{-3}F(\theta)+\theta^{-2}\xi_F-C\theta^{\nu-3})t=-\theta^{-3}(2F(\theta)+
\xi_F(-\theta)+C\theta^{\nu})t\\&\geq &-\theta^{-3}(2F(\theta)+
F^0(\theta;-\theta)+C\theta^{\nu})t\\&\geq&0,
\end{eqnarray*}
which concludes the proof. In particular, combining \eqref{eps-harmadik-bis} with \eqref{super-quadratic},  we necessarily have that  
\begin{equation}\label{nu-q}
	\nu\leq q. 	
\end{equation}

\underline{Step 4.} (Nonsmooth Cerami compactness condition for $\mathcal E_G$) Let $\{u_k\}_k\subset H^1_G(M)$ be a Cerami sequence for $\mathcal E_G$, i.e., for some $M>0$, one has  $|\mathcal E_G(u_k)|\leq M$ and $(1+\|u_k\|_{H^1})m(u_k)\to 0$ as $k\to \infty$, where $m(u)=\min\{\|\xi\|_*:\xi\in \partial \mathcal E_G(u)\}$. Our objective is to prove that, up to a subsequence, $\{u_k\}_k$ strongly converges to some element in $H^1_G(M)$. 

We first prove that $\{u_k\}_k$ is bounded in $L^\nu(M)$. For every $k\in \mathbb N$, let $\xi_k\in \partial \mathcal E_G(u_k)$ be such that $\|\xi_k\|_*=m(u_k)$. We observe that 
$$\mathcal E_G^0(u_k;u_k)\geq \langle \xi_k,u_k\rangle \geq-\|\xi_k\|_*\|u_k\|_{H^1}\geq -(1+\|u_k\|_{H^1})m(u_k).$$
Since $(1+\|u_k\|_{H^1})m(u_k)\to 0$ as $k\to \infty$, there exists $k_0\in \mathbb N$ such that for every $k>k_0$ one has that 
$\mathcal E_G^0(u_k;u_k)\geq -1.$ Consequently,  Proposition \ref{prop-lok-lip-0}/$(v)$, inequality \eqref{Rocky-type} (which is also valid due to \eqref{eps-harmadik}) and hypothesis  $(\textbf{H})_4$ imply for every $k\in \mathbb N$ that
\begin{eqnarray*}
	2M+1&\geq& 2\mathcal E_G(u_k)-\mathcal E_G^0(u_k;u_k)\\&=&\mathcal N_\mu(u_k)-2\lambda\mathcal F(u_k)-\frac{1}{2}\langle \mathcal N'_\mu(u_k);u_k\rangle-\lambda(-\mathcal F)^0(u_k;u_k)\\&=&-\lambda \left(2\mathcal F(u_k)+\mathcal F^0(u_k;-u_k)\right)\\&\geq &-\lambda \int_M \alpha(x)\left(2F(u_k(x))+F^0(u_k(x);-u_k(x))\right){\rm d}v_g\\&\geq& \lambda C \int_M \alpha(x)|u_k(x)|^\nu{\rm d}v_g.
\end{eqnarray*}
Since $\alpha\in L^\infty(M)$ and ${\rm essinf}_{x\in M}\alpha(x)=\alpha_0>0$, the latter estimate implies that
$$2M+1\geq \lambda C \alpha_0 \|u_k\|^\nu_{L^\nu},\ \forall k\in \mathbb N.$$
Consequently, $\{u_k\}_k$ is bounded in $L^\nu(M)$. 

Now, we prove that $\{u_k\}_k$ is bounded in $H_G^1(M)$.   By \eqref{eps-harmadik-bis}, for every small $\epsilon>0$ there exists $\tilde C_\epsilon>0$ such that  for every $k\in \mathbb N$,
$$M\geq \mathcal E_G(u_k)=\frac{1}{2}\mathcal N_\mu(u_k)-\lambda\mathcal F(u_k)\geq \left(\frac{c_\mu}{2}-\epsilon\lambda\|\alpha\|_{L^\infty}\right)   \|u_k\|_{H^1}^2-\lambda \tilde C_\epsilon \|\alpha\|_{L^\infty}\|u_k\|_{L^q}^q. $$
In particular, if $\frac{c_\mu}{4}>\epsilon\lambda\|\alpha\|_{L^\infty}$, then there exists $M_\epsilon>0$ and $\overline C_\epsilon>0$ such that
\begin{equation}\label{kicsi}
	\|u_k\|_{H^1}^2\leq M_\epsilon+\overline C_\epsilon\|u_k\|_{L^q}^q,\ \forall k\in \mathbb N.
\end{equation}

On account of \eqref{nu-q}, we distinguish two cases: 

a) $\nu=q$. Since $\{u_k\}_k$ is bounded in $L^\nu(M)$ and $\nu=q$, by \eqref{kicsi} we also have that  $\{u_k\}_k$ is bounded in $H_G^1(M).$

b) $\nu<q$. Let $\eta\in (0,1)$ be such that $\frac{1}{q}=\frac{1-\eta}{\nu}+\frac{\eta}{2^*}$. By \eqref{kicsi} and a standard interpolation inequality we have that
\begin{eqnarray}\label{utolsok-kozott}
\nonumber	\|u_k\|_{H^1}^2&\leq& M_\epsilon+\overline C_\epsilon\|u_k\|_{L^q}^q\leq M_\epsilon+\overline C_\epsilon\|u_k\|_{L^\nu}^{(1-\eta)q} \|u_k\|_{L^{2^*}}^{\eta q}\\&\leq& M_\epsilon+\overline C_\epsilon(K_q^\pm)^{\eta q}\|u_k\|_{L^\nu}^{(1-\eta)q} \|u_k\|_{H^1}^{\eta q},
\end{eqnarray}
where 
$K_q^\pm>0$ are the embedding constants from \eqref{Sobolev-CH} and 
\eqref{Sobolev-Ric}, respectively. Since $q<2+\frac{4}{n}$, we have that $\nu>2>\frac{n(q-2)}{2}$. We observe that $\nu>\frac{n(q-2)}{2}$ together with $\frac{1}{q}=\frac{1-\eta}{\nu}+\frac{\eta}{2^*}$ is equivalent to $\eta q<2.$ The latter inequality and \eqref{utolsok-kozott} imply that $\{u_k\}_k$ is bounded in $H_G^1(M).$ 

Now, we can proceed as in Step 4, see \S \ref{section-4}; in this way we conclude that $\{u_k\}_k$ strongly converges (up to a subsequence) to some element in $H^1_G(M)$.   

\underline{Step 5.} (Existence/multiplicity of critical points for $\mathcal E_G$) Under the assumptions of the theorem, one can prove as above that $\mathcal E_G$ has the mountain pass geometry. By Step 4 and on account of the mountain pass theorem for locally Lipschitz functions, see e.g.\ Kourogenis and Papageorgiou \cite{KP}, we conclude the existence of a non-zero critical point  for $\mathcal E_G$. When $F$ is even, we may apply the nonsmooth fountain theorem involving the Cerami compactness condition, see e.g. Krist\'aly \cite{Kristaly-SVAN}, guaranteeing the existence of a  sequence of critical points for the functional $\mathcal E_G$. All these points are $G$-invariant solutions to the differential inclusion \eqref{problem-1}, which concludes the proof. \hfill $\square$  

%
%


\begin{thebibliography}{99}
		
			\bibitem{Aubin} 
				T. Aubin, Probl\`emes isop\'erim\'etriques et espaces de Sobolev. \textit{J. Differential Geometry} 11 (1976), no. 4, 573--598.
				
				\bibitem{Berchio} E. Berchio, A.  Ferrero, G. Grillo,  Stability and qualitative properties of radial solutions of the Lane-Emden-Fowler equation on Riemannian models.\textit{ J. Math. Pures Appl.} (9) 102 (2014), no. 1, 1--35.
		
	\bibitem{Balogh-Kristaly-Annalen} Z.M.	Balogh, A.  Krist\'aly, Sharp isoperimetric and Sobolev inequalities in spaces with nonnegative Ricci curvature. \textit{Math. Annalen}, in press. DOI: 10.1007/s00208-022-02380-1
	
	\bibitem{Bonanno-etal} G. Bonanno, G. D'Agu\`i, P. Winkert, Sturm-Liouville equations involving discontinuous nonlinearities. \textit{Minimax Theory Appl.} 1 (2016), no. 1, 125--143. 
	
	\bibitem{Bonanno-Bisci} G. Bonanno, G. Molica Bisci,  V. R\u adulescu, Nonlinear elliptic problems on Riemannian manifolds and applications to Emden-Fowler type equations. \textit{Manuscripta Math.} 142 (2013), no. 1-2, 157--185.
	
	\bibitem{Candito-Livrea} P. Candito, R. Livrea,  \textit{Nonlinear difference equations with discontinuous right-hand side.} Differential and difference equations with applications, 331--339, Springer Proc. Math. Stat., 47, Springer, New York, 2013.
	
	\bibitem{Carl-Le} S. Carl, V.K. Le, 
	\textit{Multi-valued variational inequalities and inclusions}.
	Springer Monographs in Mathematics. Springer, Cham, 2021. 
	
	\bibitem{Carl-Le-2} S. Carl, V.K. Le, Extremal solutions of multi-valued variational inequalities in plane exterior domains. \textit{J. Differential Equations} 267 (2019), no. 8, 4863--4889.
	
	\bibitem{CLM} S. Carl, V.K. Le, D. Motreanu,  \textit{Nonsmooth variational problems and their inequalities}. Comparison principles and applications. Springer Monographs in Mathematics. Springer, New York, 2007.
	
	\bibitem{Chang} K.C. Chang,  Variational methods for nondifferentiable functionals and their applications to partial differential equations. \textit{J. Math. Anal. Appl.} 80 (1981), 102--129. 
		
	\bibitem{Clarke}	F.H. Clarke, \textit{Optimization and Non-smooth Analysis}, John Wiley \& Sons, New York, 1983.
	
	\bibitem{CKV} N. Costea, A.  Krist\'aly, C. Varga, \textit{Variational and Monotonicity Methods in Nonsmooth Analysis}, Frontiers in Mathematics, Birkh\"auser/Springer, 2021. 
	
	\bibitem{Coulhon-SC} T. Coulhon, L. Saloff-Coste,  Isop\'erim\'etrie pour les groupes et les vari\'et\'es. \textit{Rev. Mat. Iberoamericana} 9(2) (1993), 293--314. 
	
	\bibitem{Dipierro}  L. D'Ambrosio, S. Dipierro, Hardy inequalities on Riemannian manifolds and applications. \textit{Ann. Inst. Henri Poincar\'e, Anal. Non Lin\'eaire} 31 (3) (2014) 449--475. 
	
	\bibitem{FKM} C. Farkas, A. Krist\'aly, \'A. Mester,  Compact Sobolev embeddings on non-compact manifolds via orbit expansions of isometry groups. \textit{Calc. Var. Partial Differential Equations} 60 (2021), no. 4, Paper No. 128, 31 pp.
	
	\bibitem{Gasinski-Papa} L. Gasi\'nski, N.S. Papageorgiou, \textit{ Nonsmooth critical point theory and nonlinear boundary value problems}. Series in Mathematical Analysis and Applications, 8. Chapman \& Hall/CRC, Boca Raton, FL, 2005.
		
			\bibitem{Hebey} E. Hebey, \textit{Nonlinear analysis on manifolds: Sobolev spaces and inequalities}. Courant Lecture Notes in Mathematics, 5.
			New York University, Courant Institute of Mathematical Sciences, New
				York; American Mathematical Society, Providence, RI, 1999.	
				
				\bibitem{Jaber} H. Jaber, 
				Hardy-Sobolev equations on compact Riemannian manifolds.
			\textit{	Nonlinear Anal.} 103 (2014), 39--54. 
			
		\bibitem{KP}	N.-C. Kourogenis, N.-S. Papageorgiou, Nonsmooth critical point theory and nonlinear elliptic equations at
			resonance. \textit{Kodai Math. J.} 23 (2000) 108--135.
				
			\bibitem{KM}W.	Krawcewicz, W. Marzantowicz,  Some remarks on the Lusternik-Schnirelman method for nondifferentiable functionals invariant with respect to a finite group action. \textit{Rocky Mountain J. Math.} 20 (1990), 1041--1049.
				
				
			\bibitem{Kristaly-JMPA} A.  Krist\'aly, Sharp uncertainty principles on Riemannian manifolds: the influence of curvature. \textit{J. Math. Pures Appl.} (9) 119 (2018), 326--346.
			
				\bibitem{Kristaly-JFA} A.  Krist\'aly, New geometric aspects of Moser-Trudinger inequalities on Riemannian manifolds: the non-compact case. \textit{J. Funct. Anal.} 276 (2019), no. 8, 2359--2396.
			
			\bibitem{Kristaly-SVAN}  A.  Krist\'aly, Multiplicity results for an eigenvalue problem for hemivariational inequalities in strip-like domains.\textit{ Set-Valued Anal.} 13 (2005), no. 1, 85--103. 
			
				
				
			\bibitem{KMM} A.	Krist\'aly, \'A. Mester,  I.I. Mezei,   Anisotropic symmetrization and Sobolev inequalities on Finsler manifolds with nonnegative Ricci curvature. Preprint, 2022. Link:  https://arxiv.org/abs/2107.00512v1 
			
			\bibitem{KMV} A. Krist\'aly, V. Motreanu, C. Varga, A minimax principle with general Palais-Smale conditions. \textit{Comm. Appl. Anal.} 9(2005) (2), 285--299.
			
			\bibitem{KV} A. Krist\'aly, C. Varga, Variational-hemivariational inequalities on unbounded domains.  \textit{Stud. Univ. Babe\c s-Bolyai Math.} 55 (2010), no. 2, 3--87.
				
				
				\bibitem{Lions}  P.-L. Lions. Sym\'etrie et compacit\'e dans les espaces de Sobolev. \textit{J. Funct. Anal.} 49(3) (1982), 315--334.
				
					\bibitem{Hanne-Varga} H. Lisei, C.  Varga,  A multiplicity result for a class of elliptic problems on a compact Riemannian manifold. \textit{J. Optim. Theory Appl.} 167 (2015), no. 3, 912--927. 
				
			
				
				\bibitem{Liu-Liu} C. Liu, Y.  Liu,  The existence of ground state solution to elliptic equation with exponential growth on complete noncompact Riemannian manifold. \textit{J. Inequal. Appl.} 2020, Paper No. 74, 21 pp.
				
					\bibitem{Liu-Liu-Motreanu} Y. Liu, Z. Liu, D. Motreanu,  Differential inclusion problems with convolution and discontinuous nonlinearities. \textit{Evol. Equ. Control Theory} 9 (2020), no. 4, 1057--1071.
					
						\bibitem{olaszok-1} Z. Liu, R. Livrea, D.  Motreanu, S. Zeng, Variational differential inclusions without ellipticity condition. \textit{Electron. J. Qual. Theory Differ. Equ.} 2020, Paper No. 43, 17 pp.
						
					\bibitem{MOS}	S. Mig\'orski, A. Ochal, M. Sofonea, \textit{Nonlinear inclusions and hemivariational inequalities}. Models and analysis of contact problems. Advances in Mechanics and Mathematics, 26. Springer, New York, 2013.
				
			\bibitem{BPucci}	G. Molica Bisci, P. Pucci, \textit{Nonlinear Problems with Lack of Compactness},  De Gruyter Series in Nonlinear Analysis and Applications, 2021.
				
				\bibitem{Molica-0} G. Molica Bisci, S. Secchi,  Elliptic problems on complete non-compact Riemannian manifolds with asymptotically non-negative Ricci curvature. \textit{Nonlinear Anal.} 177 (2018), part B, 637--672.
				
				\bibitem{Molica} G. Molica Bisci, L. Vilasi,  Isometry-invariant solutions to a critical problem on non-compact Riemannian manifolds. \textit{J. Differential Equations} 269 (2020), no. 6, 5491--5519.
				
				\bibitem{Motreanu-Pan}  D. Motreanu, P.D. Panagiotopoulos,  \textit{Minimax theorems and qualitative properties of the solutions of hemivariational inequalities}. Nonconvex Optimization and its Applications, 29. Kluwer Academic Publishers, Dordrecht, 1999.
				
			\bibitem{Munn} M.	Munn, Volume growth and the topology of manifolds with nonnegative Ricci curvature. \textit{J. Geom. Anal.} 20(3) (2010), 723--750.
				
				
				\bibitem{Palais} R.S. Palais, The principle of symmetric criticality. \textit{Comm. Math. Phys.} 69 (1979), 19--30.
		
		 \bibitem{Panagiotopoulos} P.D. Panagiotopoulos, \textit{Hemivariational inequalities}. Applications in mechanics and engineering. Springer-Verlag, Berlin, 1993.
		
		\bibitem{ST-1} L. Skrzypczak, C. Tintarev, A geometric criterion for compactness of invariant subspaces. \textit{Arch. Math.} (Basel)
		101(3) (2013), 259--268.
		
		\bibitem{ST-2} L. Skrzypczak, C. Tintarev, On compact subsets of Sobolev spaces on manifolds.  \textit{Trans. Amer. Math. Soc.} 374 (2021), no. 5, 3761--3777.
		
		
		\bibitem{Talenti} G. Talenti, Best constant in Sobolev inequality.\textit{ Ann. Mat. Pura Appl.} 110 (1976), 353--372.
		
		\bibitem{Varga} C. Varga,  Existence and infinitely many solutions for an abstract class of hemivariational inequalities. \textit{J. Inequal. Appl.} 2005, no. 2, 89--105.
		
		
		
	
		

	

	
	
		
%
%
%
%
%
%
%
%
%
%
%
%
%
%
%
%
%
%
%
%
%
%
%
%
%
%
%
%
%
%
%
%
%
%
%
%
%
%
%
%
%
%
%
%
%
%
%
%
%
%
%
%
%
%
%
%
%
%
%
%
%
		
	\end{thebibliography}
\end{document}